\magnification=\magstep1
\input amstex
\documentstyle{amsppt}

\NoBlackBoxes
\loadbold
\def\mbox#1{\leavevmode\hbox{#1}}
\input btxmac
\input xy
\xyoption{v2}


\def\Index{\operatorname{Index}}
\def\Ind{\operatorname{Ind}}
\def\Ker{\operatorname{Ker}}
\def\Coker{\operatorname{Coker}}
\def\Range{\operatorname{Range}}

\def\Lim{\lim_{t \to \infty}}
\def\Prop{\operatorname{Prop}}
\def\supp{\operatorname{supp}}

\def\spinc{$\text{spin}^c$ }
\def\sym{\operatorname{sym}}
\def\Ell{\operatorname{Ell}}

\def\End{\operatorname{End}}
\def\Domain{\operatorname{Domain}}
\def\D#1{-\sqrt{-1}\ \frac{\partial}{\partial x_{#1}}}
\def\Dx#1{\frac{\partial}{\partial x_{#1}}}

\def\AM#1#2#3{\{#1_t\}_{t \in [1,\infty)} : #2 \to #3}
\def\amn#1#2#3{\{#1_t\} : #2 \to #3}
\def\i{\sqrt{-1}}
\def\im{\lrcorner}
\def\gtimes{\hat{\otimes}}
\def\dolbeaut{\partial\!\!\!/}
\def\bsigma{\boldsymbol\sigma}

\def\E{\Cal{E}}
\def\b{\Bbb}

\def\Z{\Bbb Z}
\def\Q{\Bbb Q}
\def\R{\Bbb R}
\def\C{\Bbb C}
\def\H{\Cal{H}}
\def\B{\Cal{B}}
\def\KK{\Cal{K}}
\def\L{\Cal{L}}
\def\S{\Cal{S}}


\def\lcl{[\![}
\def\rcl{]\!]}
   


\topmatter
\title Asymptotic Morphisms and Elliptic Operators over $C^*$-algebras \endtitle
\rightheadtext{Asymptotic Morphisms and Elliptic Operators}
\author Jody Trout \endauthor
\address
Department of Mathematics,
6188 Bradley Hall,
Hanover, NH 03755 \endaddress
\email Jody.Trout\@Dartmouth.EDU \endemail

\abstract
This paper provides an $E$-theoretic proof of an exact form, due to E. Troitsky, of the Mischenko-Fomenko Index Theorem for elliptic 
pseudodifferential operators over a unital $C^*$-algebra. The main ingredients in the proof are the use of asymptotic morphisms of 
Connes and Higson, vector bundle modification, a Baum-Douglas-type group, and a $KK$-argument of Kasparov.
\endabstract

\endtopmatter

\document

\head {\bf 1. Introduction} \endhead

Let $A$ be a  $C^*$-algebra with unit and $M$ be a smooth closed compact manifold. Mishchenko and Fomenko \cite{MF80} consider an
elliptic pseudodifferential $A$-operator  \linebreak $P : C^\infty(E_1) \to C^\infty(E_2)$ acting between the spaces
of smooth sections of smooth vector $A$-bundles $E_i \to M$, whose fibers are finite projective modules over $A$. The analytic index
of $P$ is the $K$-theory-valued Fredholm index $\Index_a(P) \in K_0(A)$. If the kernel and cokernel of $P$ are finite
projective $A$-modules, then it follows $\Index_a(P) = [\Ker(P)] - [\Coker(P)] \in K_0(A).$ The topological index of $P$ is also an element
$\Index_t(P) \in K_0(A)$ defined by a Gysin map construction which embeds $M$ into Euclidean space and invokes Bott Periodicity. This uses the
principal symbol
$\sigma(P) : \pi^*(E_1) \to \pi^*(E_2)$ which defines an element $[\sigma(P)] \in K^0_A(T^*M) \cong K_0(C_0(T^*M) \otimes A)$, the topological $K$-theory of vector
$A$-bundles on the cotangent bundle $\pi : T^*M \to M.$ 

Our goal is to prove that these two $K$-theory classes are actually the same
$$\Index_a(P) = \Index_t(P) \in K_0(A),$$
not in $K_0(A) \otimes \Q$ (which kills torsion) as originally proved by Mishchenko and Fomenko.
If $A = \C$, this is the classical Atiyah-Singer Index Theorem \cite{AS68}. This theorem has also been obtained by E. Troitsky
\cite{Troi96,Troi88,Troi93a,Troi93}. A complete proof is given in \cite{Troi96}, where he uses a generalization of the axiomatic method of Atiyah and Singer. In
this paper, we prove this index theorem using new $E$-theoretic asymptotic morphism techniques, which should generalize well to equivariant and graded versions of
the index theorem.

The hardest part in proving these index theorems is contained in showing that the analytic index is preserved with respect to changing the underlying base manifold
from $M$ to $S$ where $S \to M$ is a smooth compact fiber bundle as contained in the multiplicative axiom B.3
\cite{AS68}. By the Thom isomorphism $K^0_A(T^*M) @>{\cong}>> K^0_A(T^*S)$, it follows that the topological index is well-behaved with respect to this operation.
However, the analytic index is far more delicate, even in the classical setting. Our technique makes consistent use of this ``vector bundle modification''
construction (Section 2), but bypasses the difficult  calculation for the analytic index by appealing to the asymptotic morphisms of the Connes-Higson $E$-theory
\cite{CH89}. 

Specifically, we use an asymptotic morphism (constructed in Appendix A)
$$\AM {\Phi} {C_0(T^*M) \otimes A} {\Cal K(L^2M) \otimes A}$$
which is naturally associated to $M$ and $A$ up to equivalence.
If $A = \C$, this asymptotic morphism is essentially the same as the one used by Higson in his proof of the index theorem for
classical first-order differential operators \cite{Hig93}. In Section 3, we prove that the induced map on $K$-theory, 
$$\Phi_* : K^0_A(T^*M) \to  K_0(A),$$
is precisely the topological index $\Index_t(P) = \Phi_*([\sigma(P)])$. This is done by showing that the topological index and 
this ``morphism'' index $\Index_m = \Phi_*$ induce the same group isomorphism
$$\Index_t = \Index_m : \Ell(A) @>{\cong}>> K_0(A)$$
where $\Ell(A)$ is a $K$-homological Baum-Douglas-type group \cite{BD82} which incorporates vector bundle modification. 
(In Appendix B, we prove that the induced map is Bott Periodicity if $M = \R^{n}$ by generalizing Atiyah's elliptic operator proof
\cite{Ati68} and relate it to the Thom isomorphism.)

In Section 4, we develop techniques that show the above asymptotic morphism, in a sense, ``quantizes'' the principal symbol $\sigma(D)$ (considered as a
matrix-valued function on the phase space $T^*M$) of a self-adjoint first-order elliptic differential
$A$-operator $D$ on Euclidean space $M = \R^n$. The index theorem is then proven in Section 5 as follows. First, we establish it for first-order elliptic
differential $A$-operators (on arbitrary smooth closed manifolds) by adapting Higson's asymptotic morphism method for $A = \C$ mentioned above. We then use a
$KK$-theory argument of Kasparov which says that, up to homotopy, every elliptic pseudodifferential $A$-operator on a smooth
closed \spinc manifold is given by a first-order Dirac operator $D_E$ twisted by a vector $A$-bundle $E$.

The material in this paper formed a part of my Ph.D. thesis \cite{Trou95} at the Pennsylvania State
University. I want to thank my advisors Paul Baum and Nigel Higson for their
great help and encouragement. I also want to thank Guennadi Kasparov and the referee for their
helpful suggestions. 

\head {\bf 2. Vector Bundle Modification} \endhead

For a smooth closed Riemannian manifold $M$, let $L^2(M)$ denote the Hilbert space of
square-integrable functions on $M$. If $E @>p>> M$ is a complex Hermitian bundle,
let $L^2(M,E)$ (or $L^2(E)$ if there is no confusion) denote the Hilbert space of square-integrable sections of $E$. Let $C_\tau(E)$
denote the $C^*$-algebra of bundle endomorphisms of the pull-back bundle $\b E = p^*\Lambda^*E \to E$
which vanish at infinity on $E$, where $\Lambda^*E$ is the exterior algebra bundle. (See Definition B.8.)

Let $F @>{p}>> M$ be a smooth Euclidean vector bundle on the manifold $M$. Let $s : M \hookrightarrow F$ denote the
canonical embedding of $M$ into $F$ as the zero section. It follows that
$F$ is the normal bundle of this embedding. Let $s_* : TM \hookrightarrow
TF$ denote the induced embedding. The normal bundle of this embedding is
$TF$ and is just the pull-back to $TM$ of  $F \oplus F$ \cite{AS68,LM89}, that is, $TF = \pi_M^*(F \oplus F)$, where
$\pi_M : T^*M \to M$ denotes the
projection of the cotangent bundle. This bundle $TF \to TM$ has a canonical complex structure given by
$$J = \pmatrix 0 & -I \\ I & 0 \endpmatrix.$$ By using the given metrics
on $M$ and $F$, we identify $\pi : E = T^*F \to T^*M$ as this complex
bundle. We denote $C_\tau(T^*F) = C_\tau(E)$ as above. Let $F_{\C} = F \otimes_{\R} \C = F \oplus iF$ denote the
complexification of the bundle $F$. Let $\b F = p^*\Lambda^*F_{\C}$.
(Compare Definition B.8 in the appendix.) Let $\pi_F : T^*F \to F$ denote the cotangent
bundle of $F$ considered as a manifold.

\proclaim{Lemma 2.1} $\b E = \pi^*(\Lambda^*E) \cong \pi_F^*(\b F)$. \endproclaim

\demo{Proof} By definition, $E \cong \pi_M^*(F_{\C})$ as a complex vector bundle
and the following diagram $$\CD
T^*F @>{\pi}>> T^*M \\
@V{\pi_F}VV   @V{\pi_M}VV \\
F @>p>> M 
\endCD$$
commutes. Therefore, $\b E \cong \pi^*\Lambda^*\pi_M^*F_{\C} \cong
\pi_F^* \Lambda^*p^*F_{\C} = \pi_F^* \b F$. \qed \enddemo

Choose a complex vector bundle $\b G \to F$ such that $\b F \oplus \b G
\cong F \times \C^n$ is trivial \cite{Mil58}. It follows by the previous
lemma that $\b E \oplus \pi_F^*\b G \cong T^*F \times \C^n$. Therefore,
we have the following isometry of Hilbert spaces
$$V : L^2(F,\b F) \hookrightarrow L^2(F, \b F \oplus \b G) \cong
L^2(F)^{n}.$$
This induces the inclusion 
$$Ad(V) : \KK(L^2(F, \b F)) \hookrightarrow \KK(L^2(F)^{n}) =
M_n(\KK(L^2F)),$$
by the mapping $K \mapsto VKV^*$. We also have the inclusion of $C^*$-algebras
$$J : C_\tau(T^*F) \hookrightarrow C_0(T^*F, \End \C^n) \cong
M_n(C_0(T^*F)).$$

\proclaim{Lemma 2.2} {\rm 1.)} $J_* : K_0(C_\tau(T^*F)) \to K_0(C_0(T^*F))$ is an isomorphism.

\rm{2.)} ${Ad(V)_*} = id : K_0(\KK(L^2(F,\b F))) \to K_0(\KK(L^2F))$
\endproclaim

\demo{Proof} 1.) If $\b F$ is trivial, the result holds. In general, apply a Mayer-Vietoris 
argument and  the Five Lemma.

2.) Recall that any two isometries of a separable Hilbert space are
connected by a strongly continuous path of isometries. Hence, $V$ is
homotopic to a unitary isomorphism $U$. It follows that $Ad(V)$ is
homotopic to $Ad(U)$. Thus, $Ad(V)_* = Ad(U)_* = id$ since it maps rank
one projections to rank one projections. \qed \enddemo

Referring to Appendix A, let $$\{\Phi_t^M\} : C_0(T^*M) \to \KK(L^2M)$$ be the index asymptotic morphism 
for $M$. The following lemma allows us to ``twist'' this asymptotic
morphism with a complex Hermitian bundle $H \to M$. This is the ``vector bundle modification''
construction for the index asymptotic morphism. 

Let $C_H(T^*M)$ denote the $C^*$-algebra of endomorphisms of the vector bundle $\pi_M^*H \to T^*M$
vanishing at infinity in the operator norm induced by the metrics on $H$ and $M$.

\proclaim{Lemma 2.3} Let $H \to M$ be a smooth Hermitian bundle on $M$. There is an asymptotic morphism
$$\amn {\Phi^H} {C_H(T^*M)} {\KK(L^2(M, H))}$$
such that if $\alpha \in C_H(T^*M)$ has support in $T^*U$, where $H|_U \cong U \times \C^n$ on the open subset $U \subset M$, then
$$\lim_{t \to \infty} \|\Phi^H_t(f) - M_n(\Phi^M_t(f))\| = 0$$
where we identify $C_H(T^*U) \cong C_0(T^*U, \End \C^n) \cong M_n(C_0(T^*U))$.
\endproclaim

\demo{Proof} If $H \cong M \times \C^n$ is trivial, then $C_H(T^*M) \cong M_n(C_0(T^*M))$
and $L^2(M,H) \cong L^2(M)^{n}$. Define $\{\Phi^H_t\}$ to be the $n \times n$
matrix extension of $\{\Phi^M_t\}$. Now use a partition of unity $\{\rho_j^2\}$ subordinate
to an open cover $\{U_j\}_1^m$ over which $H$ trivializes $H |_{U_j} \cong U_j \times \C^n$
and a gluing argument.
\qed \enddemo

Considering $E = T^*F$ as a manifold and $H = \b F \to F$, there is, by tensoring with
the identity $id_A : A \to A$,
an asymptotic morphism
$$\amn {\Phi^{\b F, A}} {C_\tau(T^*F) \otimes A} {\KK(L^2(F, \b F)) \otimes A}.$$ We want to use
this asymptotic morphism to relate the $A$-index asymptotic morphisms
associated to the {\it manifolds} $M$ and $F$,
$$\aligned
&\amn {\Phi^{F,A}} {C_0(T^*F) \otimes A} {\KK(L^2F) \otimes A} \\
&\amn {\Phi^{M,A}} {C_0(T^*M) \otimes A} {\KK(L^2M) \otimes A}
\endaligned$$
and the Thom homomorphism (Proposition B.13)
$$\Psi = \Psi^E : C_0(\R) \otimes C_0(T^*M) \to C_0(\R) \otimes C_\tau(T^*F)$$
associated to $E = T^*F \to T^*M$. To this end, we need to construct an
elliptic operator of index one acting along the fibers $F_m$ of $p : F \to M$. 

\demo{\bf Definition 2.4} For each point $m \in M$ and $t > 0$, let 
$$B^m_t : \S(F_m, \Lambda^*F_m) \to \S(F_m, \Lambda^*F_m)$$
denote the operator constructed in Definition B.3. For each $m$, the kernel $\Ker(B^m_t)$ is one-dimensional
and spanned by the $0$-form $v \mapsto e^{-t\|v\|^2}$ (Theorem B.4).
Since $B^m_t$
is $O(n)$-equivariant, the collection
$$ \b B_t = \{B^m_t : m \in M\}$$
defines a smooth family of operators acting along the fibers of $F$.

\enddemo

\demo{\bf Definition 2.5} Define, for each $t > 0$, the map
$$\alpha_t : C_c(M) \to \Ker(\b B_t)$$ by the following formula
$$\alpha_t(f)(v) = \big(\frac{2t}{\pi}\big)^{n/4} f(p(v)) e^{-t\|v\|^2}, \quad v \in
F.$$ \enddemo

The proof of the following lemma is easy.

\proclaim{Lemma 2.6} The collection $\{\alpha_t\}$ induces a continuous family of
isometries
$$\{\alpha_t\} : L^2(M) \to L^2(F, \b F)$$
which is an isomorphism onto $\Ker (\b B_t)$.
\endproclaim

\proclaim{Corollary 2.7} There is a continuous family of injective $*$-homomorphisms
$$ \{Ad(\alpha_t)\} : \KK(L^2M) \to \KK(L^2(F, \b F)).$$
\endproclaim

Define $\{\beta_t\} : C_0(\R) \otimes \KK(L^2M) \to C_0(\R) \otimes \KK(L^2(F, \b F))$ to
be the suspension $\{ 1 \otimes Ad(\alpha_t)\}$ of $\{Ad(\alpha_t)\}$ .

\proclaim{Lemma 2.8} The induced map $\beta_* = id$ on $K$-theory. \endproclaim

\demo{Proof} $\{\beta_t\}$ is homotopic to an isomorphism since
$\{Ad(\alpha_t)\}$ is homotopic to an isomorphism. \qed \enddemo

We come to the main result of this section. Consider the following
suspended diagram of $C^*$-algebras and asymptotic morphisms:

$$\CD
{C_0(\R) \otimes M_n(C_0(T^*F))} @>{\{1 \otimes \Phi^F_t\}}>> {C_0(\R)
\otimes M_n(\KK(L^2F))}\\
@AA{1 \otimes J}A                       @AA{1 \otimes Ad(V)}A \\
{C_0(\R) \otimes C_\tau(T^*F)}     @>{\{1 \otimes \Phi^{\b F}_t\}}>> {C_0(\R)
\otimes \KK(L^2(F, \b F))} \\
@AA{\Psi}A                           @AA{\{\beta_t\}}A \\
{C_0(\R) \otimes C_0(T^*M)}  @>{\{1 \otimes \Phi^M_t\}}>> {C_0(\R)
\otimes \KK(L^2M)}
\endCD$$

\proclaim{Theorem 2.9} The diagram above commutes up to homotopy.
\endproclaim

\demo{Proof} 1.) If $M = \{pt\}$, the theorem is true by Corollary B.21.

2.) If $F = M \times \C^n$, then we have that 
$$C_\tau(T^*F) = C_0(T^*M \times T^*\R^n, \End \Lambda^*\C^n) \cong
M_n(C_0(T^*\R^n)) \otimes C_0(T^*M).$$ Thus, the diagram reduces to the
following:
$$\diagram
C_0(\R) \otimes M_n C_0(T^*\R^n) \otimes C_0(T^*M) \rto^{1 \otimes
\Phi_t \otimes \Phi^M_t} & C_0(\R) \otimes M_n \KK(L^2\R^n) \otimes
\KK(L^2M) \\
C_0(\R) \otimes C_0(T^*M) \uto_{1 \otimes \Psi} \rto^{1 \otimes
\Phi^M_t} & C_0(\R) \otimes \KK(L^2M) \uto_{\beta_t \otimes 1}
\enddiagram$$
The homotopy is  given by the formula
$$f \otimes g \mapsto f(\epsilon x + s^{-1}\bar{B}_t) \otimes \Phi^M_t(g),\
0 \leq s \le 1,$$ and is the tensor product of the homotopy in Lemma B.20
with the index asymptotic morphism for the manifold $M$.

3.) In general, the homotopy in the previous case patches together via a
partition of unity from $M$ to form the desired homotopy since it is
diffeomorphism invariant by Theorem B.4 (6) and Corollary A.12. \qed \enddemo

\proclaim{Corollary 2.10} For any $C^*$-algebra $A$, the following diagram commutes:
$$\CD
{K_0(C_0(T^*F) \otimes A)}   @>{\Phi^{F,A}_*}>>  K_0(A) \\
@AA{\{\Psi \otimes id_A\}_*}A                   @AA=A  \\
{K_0(C_0(T^*M) \otimes A)}   @>{\Phi^{M,A}_*}>>  K_0(A) 
\endCD$$
\endproclaim

\demo{Proof} Tensor the diagram prior to Theorem 2.9 with $id_A : A
\to A$. The resulting diagram commutes up to homotopy by tensoring the
homotopy from Theorem 2.9 with $id_A$. Now take the induced maps and invoke
homotopy invariance. \qed \enddemo

\head {\bf 3. The Topological and Morphism Indices} \endhead

Let $A$ be a $C^*$-algebra with unit. In this section, we assemble (the symbols of) all elliptic pseudodifferential $A$-operators on all manifolds $M$ into
an abelian group $\Ell(A)$ using $K$-homological ideas of Baum and Douglas \cite{BD82}. 
The topological index and ``morphism'' index will define two group homomorphisms
$$\Ind_t, \ \Ind_m : \Ell(A) \to K_0(A).$$
By checking examples for $M = \{pt\}$ 
we will show that they are both group isomorphisms. Bott Periodicity will then show
that they are the same $\Ind_m = \Ind_t$.

Let $M$ be a smooth compact manifold. Recall that
a {\it vector $A$-bundle} $E \to M$ is a locally trivial fiber bundle whose fibers $E_p$ for $p \in M$ are given by a finite
projective $A$-module $P$ \cite{Kar78,MF80}. (The structure group is the automorphism group $Aut_A(P)$ of $P$.)
Denote by $K^0_A(M)$ the Grothendieck group of all (isomorphism classes of) vector $A$-bundles on $M$ under direct sum. 
(If $A = \C$, this is the ordinary topological $K$-theory of $M$.) If $M$ is only locally
compact, then  we can identify $K^0_A(M) =_{def} \Ker\{K^0_A(M^+) \to K^0_A(\{\infty\}) = K_0(A)\}$ as the abelian group generated by 
vector $A$-bundle homomorphisms $\sigma : E \to F$ with compact support under direct sum, where 
$\supp(\sigma) = \{p \in M \ | \ \sigma_p : E_p \to F_p \text{ is not a module isomorphism} \}.$
There are natural notions of isomorphism and pull-backs as in the classical case $A = \C$.
Since $M$ is a smooth manifold, we may take $E, F$, and $\sigma$ to also be smooth. There is the
Mingo-Serre-Swan isomorphism $K^0_A(M) \cong K_0(C_0(M)
\otimes A)$, which is induced by taking sections as in the classical case \cite{Min82,Swa62}.

A {\it Hermitian $A$-metric} on $E$ is a smooth choice $\langle
\ ,\ \rangle_p : E_p \times E_p \to A$ of Hilbert $A$-module structures on the
fibers of $E$. Every $A$-bundle has a Hermitian $A$-metric by a smooth partition of unity
argument and using the fact that any finite projective $A$-module $P$ has a
canonical Hilbert $A$-module structure (up to unitary isomorphism) \cite{WO93}. Moreover, any two
such Hermitian $A$-metrics are homotopic to each other via the straight line homotopy.

Let $\sigma : E \to F$ be a homomorphism of vector $A$-bundles
equipped with Hermitian $A$-metrics. There is an adjoint
homomorphism $\sigma^* : F \to E$ such that $$\langle \sigma_p(e_p),
f_p\rangle_p = \langle e_p, \sigma_p^*(f_p)\rangle_p$$ for all $e_p \in
E_p$ and $f_p \in F_p$. Furthermore, $\supp(\sigma^*) = \supp(\sigma)$
and $\sigma^*$ is well-defined up to homotopy.

Let $\pi : V \to M$ be a complex vector bundle. Let $c : V \to \End(\b V)$
denote the canonical section of the bundle $\b V = \pi^*(\Lambda^*V) = \b
V^{even} \oplus \b V^{odd} \to V$. Put a Hermitian metric on $V$. Define a bundle morphism $\lambda_V : \b
V^{even} \to \b V^{odd}$ by the formula $$\lambda_V(\omega_v) =
c(v)\omega_v = v \wedge \omega_v - v \im \omega_v,$$
It follows that $\supp(\lambda_V)$ is the zero section $M$ of $V$. (Compare Definition B.8 and
Lemma B.9.)

\demo{\bf Definition 3.1}{\rm (Sharp Product) For an $A$-bundle homomorphism $\sigma : E \to F$ on
$M$, define the $A$-bundle homomorphism on (on $T^*M$)
$$\pi^*(\sigma) \# \lambda_V : \pi^*E \otimes \b V^{even} \oplus \pi^*F
\otimes \b V^{odd} \to \pi^*F \otimes \b V^{odd} \oplus \pi^*E \otimes \b
V^{even}$$ by the formula 
$$\pi^*(\sigma) \# \lambda_V = 
\pmatrix \pi^*(\sigma) \otimes 1 & -1 \otimes \lambda_V^* \\
1 \otimes \lambda_V & \pi^*(\sigma^*) \otimes 1 \endpmatrix,$$
where $\pi^*(\sigma) : \pi^*E \to \pi^*F$ denotes the pull-back of
$\sigma$. Note that $$\supp(\pi^*(\sigma) \# \lambda_V) = \pi^{-
1}(\supp(\sigma)) \cap \supp(\lambda_V) = \supp(\sigma)$$ is a compact subset of $V$,
where we identify $M$ as the zero section of $V$. \enddemo

\proclaim{Proposition 3.2} {\rm (Thom Isomorphism \cite{Troi88,Troi93})} The map $\Theta : K^0_A(M) \to K^0_A(V)$ defined by
$$[\sigma] \mapsto [\pi^*(\sigma) \# \lambda_V]$$ is an isomorphism.
\endproclaim

\demo{\bf Definition 3.3} Define $\Ell(A)$ to be the Grothendieck group of the semigroup
generated by all pairs of the form $(M,\sigma)$, where $M$ is a smooth, second-countable, 
Hausdorff manifold without boundary and $\sigma : E \to F$ is a homomorphism of smooth vector 
$A$-bundles on the cotangent bundle $T^*M$ with compact support, 
subject to the following relations:
\roster 
\item {\bf Isomorphism:} $(M,\sigma) = (M',\sigma')$ if there is a
diffeomorphism $\phi : M \to M'$ such that $\sigma \cong
\tilde{\phi}^*(\sigma')$, where $\tilde{\phi} : T^*M \to T^*M'$ denotes the
induced map. \footnote{See the discussion preceding Proposition A.7}
\item {\bf Homotopy:} $(M, \sigma_0) = (M, \sigma_1)$ if $\sigma_0$ and $\sigma_1$
are homotopic ($\sigma_0 \sim_h \sigma_1$). That is , there is a pair $(M \times [0,1],
\hat{\sigma})$ such that $\hat{\sigma}|_{T*M \times 0} \cong \sigma_0$ and $\hat{\sigma}|_{T*M \times 1} \cong
\sigma_1$.
\item {\bf Excision:} $(M,\sigma) = (U, \sigma|_{T^*U})$ if $U \subset M$
is open with $\supp(\sigma) \subset T^*U$.
\item {\bf Direct Sum - Disjoint Union:} 
$$\aligned &(M, \sigma) + (M', \sigma') = (M \amalg M', \sigma \amalg \sigma') \\
&(M, \sigma) + (M, \tau) = (M, \sigma \oplus \tau)
\endaligned$$
\item {\bf Vector Bundle Modification:} If $F \to M$ is a smooth Euclidean vector
bundle, then 
$$(M, \sigma) = F\#(M,\sigma) =_{def} (F,\pi^*(\sigma)\#\lambda_V)$$
where $V = T^*F @>\pi>> T^*M$ has the complex bundle structure from the beginning of Section 2 and
$\pi^*(\sigma)\#\lambda_V$ is the sharp product in Definition 3.1.
\endroster
\enddemo

Note that $\Ell(A)$ is abelian by the isomorphism relation. A pair $(M,\sigma)$ in $\Ell(A)$ will be called a {\it symbol pair}
and $\sigma$ will be called an {\it $A$-symbol}. If $\supp(\sigma) =
\emptyset$, then $(M,\sigma)$ will be called a {\it trivial pair}. From the
definition of $\Ell(A)$, we see that 
$(M,\sigma) \in \Ell(A) \text{ if and only if } [\sigma] \in K^0_A(T^*M).$

\proclaim{Lemma 3.4} The following are true:

\roster
\item If $(M, \sigma)$ is trivial, then $(M,\sigma) = 0 \in \Ell(A)$.
\item  If $(M,\sigma_0) = (M,\sigma_1)$ via homotopy, then $[\sigma_0] = [\sigma_1] \in K^0_A(T^*M)$. 
\item If $[\sigma] = 0 \in K^0_A(T^*M)$, then $(M,\sigma) = 0 \in
\Ell(A)$.
\endroster
\endproclaim

\proclaim{Proposition 3.5} For every symbol pair $(M,\sigma)$, there is a
pair $(\R^n, \tau)$ such that $(M,\sigma) = (\R^n,\tau) \in \Ell(A)$, for
some $n$. \endproclaim

\demo{Proof} Let $\phi : M \hookrightarrow \R^n$ be a smooth embedding
of the manifold into Euclidean space, which exists by the Whitney Embedding Theorem. By Isomorphism (1),
we may identify $M$ with its image in $\R^n$. Let $N$ be an open tubular neighborhood of $M$ in
$\R^n$. Then $M \subset N \subset \R^n$ and $N$ has the structure of an
$\R$-vector bundle $\pi : N \to M$. By Vector Bundle Modification (5), we
have that
$$(M,\sigma) = (N, \hat{\sigma})$$
where $\hat{\sigma} = p^*(\sigma) \# \lambda_{T^*N} : E \to F$ ($E$ and
$F$ are vector $A$-bundles over $T^*N$) and $p : T^*N \to T^*M$. Choose an
$A$-bundle $G$ such that $F \oplus G \cong T^*N \times A^{m}$. Then
$$\hat{\sigma} \oplus 1_G : E\oplus G \to F \oplus G \cong T^*N \times
A^{m}$$
and has compact $\supp(\hat{\sigma} \oplus 1_G) = \supp(\sigma) \subset
T^*N$. By Lemma 3.4 and Direct Sum - Disjoint Union (4), we have
$$(N,\hat{\sigma}) = (N,\hat{\sigma})+(N,1_G) = (N, \hat{\sigma} \oplus
1_G).$$
Let $H = T^*\R^n \times A^{m}$. Since $\hat{\sigma} \oplus 1_G$ is an
isomorphism off the compact set $K =\supp(\hat{\sigma})$, we can use
a clutching construction with $1_H : H \to \H$ to
obtain the $A$-symbol
$$\tau = (\hat{\sigma} \oplus 1_G) \amalg_{T^*N\backslash K} 1_H : (E
\oplus G) \amalg_{T^*N\backslash K} H \to H.$$ Since $\tau|_{T^*N} =
\hat{\sigma} \oplus 1_G$ and $\supp(\tau) \subset T^*N$, we have by
Excision (3) that
$$(M,\sigma) = (N, \hat{\sigma}) = (N,\hat{\sigma} \oplus 1_G) =
(N,\tau|_{T^*N}) = (\R^n,\tau)$$ as was desired. \qed \enddemo

Although the Gysin construction is well-known in the classical (complex bundle) setting \cite{AS68,LM89},
we will define it, since we will need to refer to the construction later.

\proclaim{Lemma 3.6} Let $g : M \hookrightarrow M'$ be a proper smooth
embedding of smooth manifolds. There is a canonical functorial
homomorphism $$g_* : K^0_A(T^*M) \to K^0_A(T^*M').$$ That is, if $h : N
\hookrightarrow Z$ is an embedding, then $(h \circ g)_* = h_* \circ g_*$.
\endproclaim

The induced map $\hat{g} : T^*M \to T^*M'$ embeds $T^*M$ as an open submanifold of $T^*M'$. Let $N(T^*M)$
be the normal bundle of this embedding (defined by pulling back the normal bundle $N(M)$ of $M$ in $M'$). It has
the structure of a smooth complex vector bundle $N(T^*M) \to T^*M$ and is an open subset of $T^*M'$. The 
map $g_*$ is defined  as the composition
$$K^0_A(T^*M) @>{\Theta}>> K^0_A(N(T^*M)) @>{i_*}>> K^0_A(T^*M')$$
where $\Theta$ is the Thom isomorphism and $i : N(T^*M) \hookrightarrow T^*M'$ is the inclusion.

\demo{\bf Remark 3.7} If $g$ is the inclusion of $M$ as an open submanifold of $N$, then
$g_*$ is the map induced by the open inclusion $T^*M \subset T^*N$, since the normal
bundle is the zero vector bundle $T^*M \times 0 = T^*M$. \enddemo

\proclaim{Proposition 3.8} If $f : M \to N$ is a smooth map of smooth manifolds,
there is a canonical homomorphism $f_! : K^0_A(T^*M) \to K^0_A(T^*N)$
depending only on the homotopy class of $f$.
\endproclaim

Let $g : M \hookrightarrow \R^n$ be any smooth embedding. Then $f \times g : M \to N \times \R^n$
is also an embedding. The {\it Gysin map} $f_!$ is defined to be the composition
$$K^0_A(T^*M) @>{(f \times g)_*}>> K^0_A(T^*N \times T^*\R^n) @>{\Theta^{-1}}>> K^0_A(T^*N)$$
and is independent of the choice of embedding $g$.

Let $f^M : M \to \{pt\}$ denote the unique map to a point.
The previous results imply that $f^M_! = \alpha_A \circ g_*$, where $g : M
\hookrightarrow \R^n$ is {\it any} embedding of $M$ into Euclidean space 
and $\alpha_A : K_0(C_0(T^*\R^{n}) \otimes A)  @>\cong>> K_0(A)$ is Bott Periodicity (Theorem B.1).

\demo{\bf Definition 3.9} {(The Topological Index)} \hfill \newline
For each $(M,\sigma) \in \Ell(A)$, define the {\bf topological index} of $(M,\sigma)$ by the formula
$$\Ind_t(M,\sigma) = f^M_!([\sigma]) \in K^0(A)$$ where $f^M_! :
K^0_A(T^*M) \to K^0_A(\{pt\}) = K_0(A)$ is the associated Gysin map.
\enddemo

\proclaim{Theorem 3.10} $\Ind_t : \Ell(A) \to K_0(A)$ induces an isomorphism of
abelian groups.
\endproclaim

\demo{Proof} First, we must show that $\Ind_t$ is well-defined, i.e., it
respects the equivalence relations in Definition 3.3 used to define the
group $\Ell(A)$.

\roster
\item {Isomorphism:} Obvious.

\item {Homotopy:} Follows from Lemma 3.4 (2). If $\sigma_0 \sim_h
\sigma_1$ are homotopic then $[\sigma_0] = [\sigma_1] \in K^0_A(T^*M)$ and so
$$\Ind_t(M,\sigma_0) = f^M_![\sigma_0] = f^M_![\sigma_1] =
\Ind_t(M,\sigma_1).$$

\item {Excision:} Suppose $U \subset M$ is open and $\supp(\sigma) \subset
T^*U$. It follows (Remark 3.7) that since $U$ is open and the inclusion $i : U
\hookrightarrow M$ is an embedding, that the normal bundle $N(T^*U) =
T^*U$ is the zero vector bundle over $T^*U$. Thus, the Thom isomorphism
$\Theta$ is the identity map on $K^0_A(T^*U)$. Let $g : M \hookrightarrow
\R^n$ be an embedding of $M$ into $\R^n$. Then $g \circ i : U
\hookrightarrow \R^n$ is an embedding. Hence we have (by Lemma 3.6) that
$$f^U_! = \Theta^{-1} \circ (g \circ i)_* = \Theta^{-1} \circ g_* \circ i_* =
f^M_! \circ i_*.$$ which implies
$$\aligned \Ind_t(U, \sigma|_{T^*U}) &= f^U_!([\sigma|_{T^*U}]) =
f^M_!(i_*([\sigma|_{T^*U}])) \\
&= f^M_!(i_* \circ i^*([\sigma]) = f^M_!([\sigma]) = \Ind_t(M,\sigma).
\endaligned$$

\item {Direct Sum-Disjoint Union:} Let $(M,\sigma)$ and $(M,\tau)$ be given.
Then we have that 
$$\aligned \Ind_t((M,\sigma) + (M,\tau)) = f^M_!([\sigma \oplus \tau])
&= f^M_!([\sigma]) + f^M_!([\tau]) \\
&= \Ind_t(M,\sigma) + \Ind_t(M,\tau).
\endaligned$$

Suppose $(N,\rho)$ is given. Since $T^*(M \amalg N) = T^*M \amalg T^*N$,
we have that $$K^0_A(T^*(M \amalg N)) = K^0_A(T^*M) \oplus
K^0_A(T^*N).$$  Let $g : M \amalg N \hookrightarrow \R^n$ be an
embedding. Consider the diagram of inclusions
$$\diagram
M \rto^-{i} \drto & M \amalg N \dto^g & N \lto_-{j} \dlto \\
                  &    \R^n           &  
\enddiagram$$
It follows by functoriality (Lemma 3.6) that
$$\aligned \Ind_t(M \amalg M, \sigma \amalg \rho) &= \Theta^{-1} \circ
g_*(i_*[\sigma] + j_*[\rho]) \\
&= f^M_!([\sigma]) + f^N_!([\rho]) = \Ind_t(M,\sigma) + \Ind_t(N,\rho).
\endaligned$$ 

\item{Vector Bundle Modification:} Let $F \to M$ be a smooth $\R$-vector
bundle. Let $s : M \hookrightarrow F$ denote the inclusion of $M$ as the
zero section. Choose an embedding  $g : F \hookrightarrow \R^n$. Then $g
\circ s : M \hookrightarrow \R^n$ is an embedding of $M$ as a closed
submanifold. By construction $s_* : K^0_A(T^*M) \to K^0_A(T^*F)$ is the
Thom Isomorphism for the complex bundle $T^*F \to T^*M$, since $F =N$ is
the normal bundle of $M$. Thus,
$$\aligned \Ind_t(F, \pi^*(\sigma) \# \lambda_{T^*F}) &= \Theta^{-1}
\circ g_* ([\pi^*(\sigma)\#\lambda_{T^*F}]) = g_*(s_*([\sigma])) \\
&= \Theta^{-1} \circ (g \circ s)_*([\sigma]) = f^M_!([\sigma]) =
\Ind_t(M,\sigma). \endaligned$$
\endroster

Therefore, $\Ind_t : \Ell(A) \to K_0(A)$ induces a well-defined homomorphism.
\footnote{Recall that
if $\phi : S \to G$ is a homomorphism from an abelian semigroup $S$ to an abelian group $G$
there is a canonical extension $\hat{\phi} : G(S) \to G$ on the Grothendieck group $G(S)$.}

Let $M =\{pt\}$ denote the unique $0$-dimensional
manifold. Then $f^M : \{pt\} \to \{pt\}$ is the identity and so, by
definition, $$f^M_! = id : K^0_A(\{pt\}) = K_0(A) \to K_0(A)$$ is the identity.
Thus, $\Ind_t$ is surjective. Now to show it is also injective.

Suppose $\Ind_t(M,\sigma) = 0$. By Proposition 3.5, we may
assume that $(M,\sigma) = (\R^n,\tau)$. By construction $$f^{\R^n}_! :
K^0_A(T^*\R^n) \to K_0(A)$$ is the (inverse) of the Thom isomorphism  for
the bundle $\C^n \cong T^*\R^n \to \{pt\}$. This implies that $[\tau] = 0 \in
K^0_A(T^*\R^n)$. By Lemma 3.4 (3), we have that $$(M,\sigma) =
(\R^n,\tau) = 0 \in \Ell(A).$$ This completes the proof that the topological
index is an isomorphism. \qed \enddemo

For any manifold $M$, let $\amn {\Phi^{M,A}} {C_0(T^*M) \otimes A}
{\KK(L^2M) \otimes A}$ denote the $A$-index asymptotic morphism of
$M$ from Appendix A. Let $$\Phi^{M,A}_* : K^0_A(T^*M) \to K_0(A)$$ denote the induced map, 
where we identify $K^0_A(T^*M) \cong K_0(C_0(T^*M) \otimes A)$ and
$K_0(\KK \otimes A) = K_0(A)$.

\demo{\bf Definition 3.11} {(The Morphism Index)} \hfill \newline
For each $(M,\sigma) \in \Ell(A)$, define the {\bf morphism index} of $(M,\sigma)$ by the formula
$$\Ind_m(M,\sigma) = \Phi^{M,A}_*([\sigma]) \in K_0(A).$$
\enddemo

\proclaim{Proposition 3.12} $\Ind_m : \Ell(A) \to K_0(A)$ induces an isomorphism of
abelian groups. \endproclaim

The proof proceeds as for the topological index.

\demo{Proof} First, we show that $\Ind_m$ respects the equivalence relations in $\Ell(A)$.

\roster
\item {Isomorphism:} Follows from diffeomorphism invariance (Corollary A.12).

\item {Homotopy:} Obvious.

\item {Excision:} Follows from the restriction property (Corollary A.11).

\item {Direct-Sum Disjoint-Union:} The relation
$$\Ind_m(M, \sigma \oplus \tau) = \Ind_m(M, \sigma) + \Ind_m(M,\tau)$$
follows since $\Phi^{M,A}_*$ is a group homomorphism. The disjoint union
relation follows from the (compatible) canonical isomorphisms 
$$\aligned C_0(T^*M \amalg T^*N) \otimes A &=(C_0(T^*M) \otimes A) \oplus (C_0(T^*N) \otimes A) \\
           \KK(L^2(M \amalg N)) \otimes A  & = (\KK(L^2M) \otimes A) \oplus (\KK(L^2N) \otimes A)
\endaligned$$
and so $\{\Phi_t^{M \amalg N, A}\} = \{\Phi_t^{M,A}\} \oplus \{\Phi_t^{N,A}\}$. 

And last but not least, the very important:

\item {Vector Bundle Modification:} Follows from Corollary 2.10 and Theorem B.22.

\endroster
Hence, $\Ind_t$ induces a well-defined group homomorphism.
If $M = \{pt\}$, then $\{\Phi_t^{M,A}\} \cong id_A : \C \otimes A \cong A \to A$ and so $\Ind_t$ is surjective.
If $M = \R^n$, then Bott Periodicity (Theorem B.7) again shows that it is injective. \qed \enddemo

The Bott Periodicity arguments in Theorem 3.10 and Proposition 3.12 show that these two isomorphisms are, in fact, the same.

\proclaim{Theorem 3.13} For any $(M,\sigma) \in \Ell(A)$, the topological and morphism indices are equal, i.e.,
$$\Ind_t(M,\sigma) = \Ind_m(M,\sigma) \in K_0(A).$$
\endproclaim

\head {\bf 4. Elliptic $A$-Operators on Euclidean Space} \endhead

Let $A$ be a $C^*$-algebra with unit. Let $D$ be a differential $A$-operator of order 
one on $\R^n$ acting on $A^k$-valued functions. Thus, $D$ has the form
$$D = \sum_{j=1}^{n} a_j(x) \Dx j + b(x),$$
where $a_j, b : \R^n \to M_k(A)$ are smooth matrix $A$-valued functions on
$\R^n$. Initially, $\Domain(D) = C_c^\infty(\R^n, A^k)$ is the pre-Hilbert $A$-module of
smooth, compactly supported vector $A$-valued functions on $\R^n$, where $A^k$ has the standard
Hilbert $A$-module structure arising from $A$.

If $D$ is formally self-adjoint on $\Domain(D)$, this implies
that the principal symbol 
$$\sigma(D)(x, \xi) = \sqrt{-1} \sum_{1}^{n} a_j(x) \xi_j \in M_k(A),$$
where $\xi = \sum_1^n\xi_j dx_j  \in
T_x^*\R^n \cong \R^n$, is a self-adjoint matrix for all $(x,\xi) \in
T^*\R^n$. $D$ is called {\it elliptic} if $\sigma(D)(x,\xi)$ is an invertible matrix for all $\xi \not = 0$.
We will also need the {\it total symbol} $\sym(D) = \sigma(D) + b$.

Consider $D$ as a densely defined unbounded symmetric operator on the
Hilbert $A$-module $\H = L^2(\R^n, A^k)$ of square-integrable vector $A$-valued
functions. It follows that $D$ has a densely defined adjoint $D^*$ and $(D \pm i)$
are injective with closed ranges. 
The operator $D$ is called {\it regular} if $(1 + D^*D)$ has dense range.
Since $D$ is symmetric, $D$ is regular if and only if $\Range(D \pm i)$ is complemented
in $\H$.  (For a review of unbounded operators on Hilbert $A$-modules, see Lance \cite{Lan95}.)

\proclaim{Theorem 4.1 (Functional Calculus)} If $D$ is a self-adjoint regular $A$-operator on the Hilbert $A$-module $\H$, 
there is a unique $*$-homomorphism
$$\phi_D : C_0(\R) \to \Cal L(\H) : f \mapsto f(D)$$
defined by sending $(x \pm i)^{-1} \mapsto (D \pm i)^{-1}$. 
\endproclaim

\definition{Definition 4.2} Let $D$ be a formally self-adjoint elliptic partial differential $A$-operator of order one on $\R^n$.
We shall call $D$ an {\it essential $A$-operator} if the closure $\bar{D}$ is self-adjoint and $(D \pm i)$ have dense range.
This implies that $\bar{D}$ is also regular. \enddefinition

Let $\phi$ be a smooth function on $\R^n$. Considered as a multiplication
operator on the Hilbert $A$-module $\H$, $M_\phi$ maps $\Domain(D)$ into itself, so
the commutator $[D, M_\phi]$ is defined on the domain of $D$. We have the
{\it symbol identity}: 
$$\i[D,M_\phi]f(x) = \sigma(D)(x, d\phi(x))f(x) = \sum_1^na_j(x)\frac{\partial \phi}{\partial x_j}f(x)$$
where $f \in \Domain(D)$. Thus, the commutator $[D, M_\phi]$ is a
pointwise skew-adjoint multiplication operator with same domain as $D$.

Using the symbol identity and the fact that if $\phi$ is a smooth and compactly supported
function then $d\phi = \sum_1^n \frac{\partial \phi}{\partial x_i} dx_i$ is a compactly supported one-form,
we obtain the following.

\proclaim{Lemma 4.3} For any smooth, compactly supported function $\phi$ on $\R^n$,
the commutator $[D,M_\phi]$ extends to a bounded adjointable $A$-operator on $\H$. In
fact,
$$\|[D,M_\phi]\| \leq \sup\{\|d\phi(x)\| \Prop(D, x) : x \in \R^n\}$$
where $\Prop(D, x) = \sup\{\|\sigma(D)(x,\xi)\| : (x,\xi) \in T^*\R^n, \|\xi\| = 1\}$.
\endproclaim

The proof of the next result is the same as in the classical case by using the generalized Rellich Lemma (Lemma 3.3 \cite{MF80}).
(The case for $A=\C$ is proved in Lemma 3.5 \cite{Gue98}.)

\proclaim{Lemma 4.4} Let $D$ be an essential $A$-operator on $\R^n$. For any 
$f \in C_0(\R)$ and $\phi \in C_0(\R^n)$, the operator $M_\phi f(\bar{D})$ is a compact $A$-operator on $\H$. \endproclaim

\proclaim{Lemma 4.5} Let $D_1$ and $D_2$ be essential $A$-operators on $\R^n$. Let $\phi \in C_0(R^n), f \in C_0(\R^n)$.  
Then we have that:
\roster
\item $\Lim \|[M_\phi, f(t^{-1}\bar{D_1})]\| = 0.$

\item If $D_2 - D_1$ is order zero, then $\Lim \|M_\phi f(t^{-1}\bar{D}_1)
- M_\phi f(t^{-1}\bar{D}_2) \| = 0.$

\item If $D_1 = D_2$ near $\supp (\phi)$, then $\Lim \|M_\phi f(t^{-
1}\bar{D}_1) -  f(t^{-1}\bar{D}_2)M_\phi\| = 0.$

\endroster

\endproclaim

\demo{Proof} 
The collection of $f \in C_0(\R)$ for which the Lemma holds is a $C^*$-
subalgebra of $C_0(\R)$. Thus, we need only check on the resolvent
functions $r_\pm (x)= (x \pm i)^{-1}$. Note that $\|r_\pm(t^{-1}\bar{D_1})\| \leq 1$
independently of $t$. 

\roster
\item From the commutator identity
$$[M_\phi, r_\pm(t^{-1}\bar{D_1})] = [M_{\phi-\phi_n}, r_\pm(t^{-1}\bar{D_1})] + [M_{\phi_n}, r_\pm(t^{-1}\bar{D_1})]$$ 
we obtain the inequality
$$\|[M_\phi,  r_\pm(t^{-1}\bar{D_1})]\| \leq 2\|\phi - \phi_n\| \| f\| + \|[M_{\phi_n}, r_\pm(t^{-1}\bar{D_1})]\|$$ 
and so we may assume that $\phi\in C^\infty_c(\R^n)$. Now, the commutator identity
$$[M_{\phi}, r_\pm(t^{-1}\bar{D_1})] = t^{-1}r_\pm(t^{-1}\bar{D_1})[D_1,M_{\phi}]r_\pm(t^{-1}\bar{D_1})$$
implies that 
$$\|[M_\phi, r_\pm(t^{-1}\bar{D_1})]\| \leq t^{-1} \|[D_1,M_{\phi}]\| \to 0$$ 
as $t \to \infty$ and the result follows.

\item Suppose $D_2 - D_1  = b(x)$, where $b$ is a smooth function. Then
we have that 
$$\aligned 
& M_\phi r_\pm(t^{-1}\bar{D}_1) - r_\pm(t^{-1}\bar{D}_2) M_\phi 
= t^{-1}r_\pm(t^{-1}\bar{D}_2)(D_2 M_\phi - M_\phi D_1)r_\pm(t^{-
1}\bar{D}_1) \\
& = t^{-1}r_\pm(t^{-1}\bar{D}_2)(M_{\phi b}+[D_2, M_\phi])r_\pm(t^{-
1}\bar{D}_1)   \to 0 \text{ as } t \to \infty
\endaligned$$  since $M_{\phi b}$ and $[D_2,M_\phi]$ are
bounded. The result follows.

\item Suppose $D_2 - D_1 = 0$ near $\supp (\phi)$. Then $M_\phi(D_2 -
D_1) = 0$ and the result follows from the previous computation. \qed 
\endroster \enddemo 

\proclaim{Theorem 4.6} Let $D$ be an essential $A$-operator on $\R^n$.
There is an asymptotic morphism $\E^D_t : \C_0(\R) \otimes C_0(\R^n) \to
\KK(\H)$
which is determined (up to asymptotic equivalence) by the maps
$$f \otimes \phi \mapsto M_\phi f(t^{-1}\bar{D}).$$
The bounded operator $f(t^{-1}\bar{D})$ on $\H$ is defined via the
functional calculus.
\endproclaim

Let $\AM {\Phi^A} {C_0(\R^{2n}) \otimes A } {\KK(L^2\R^n) \otimes A}$ be the $A$-index asymptotic
morphism defined in Appendix A.
 Recall that Banach $A$-module $L^1(\R^n, A)$ is an algebra under convolution
$$f * g(x) = \int_{\R^n} f(x-y)g(y) dy.$$
It also acts as adjointable $A$-operators via this formula on the Hilbert $A$-module $L^2(\R^n, A)$. Let $C^*(\R^n, A)$ denote
the $C^*$-algebra completion of $L^1(\R^n, A)$ in the operator norm of $\L(L^2(\R^n,A)).$ The generalized Fourier Transform
$$\hat{f}(\xi) = \frac{1}{(2\pi)^{n/2}} \int_{\R^n} f(x) e^{-ix \xi} dx$$
determines an isomorphism $\wedge : C^*(\R^n, A) \to C_0(\R^n, A) \cong C_0(\R^n) \otimes A$.
We can define a continuous family of $*$-homomorphisms 
$$C_t^A : C_0(\R^n, A) \to \L(L^2(\R^n,A)) : [C_t^A(f)]\eta = \check{f_t} * \eta.$$
The next result follows easily from Lemma A.1 since $C_t^A = C_t \otimes id_A$ with
respect to the isomorphism $C_0(\R^n, A) \cong C_0(\R^n) \otimes A$.

\proclaim{Lemma 4.7} For  any $\phi \in C_0(\R^n)$ and $f \in  C_0(\R^n, A)$, we have 
$$\lim_{t \to \infty} \|\Phi_t^A(\phi \otimes f) - M_{\phi}C_t^A(f)\| = 0$$
\endproclaim

Extend $\{\Phi_t^A\}$ to $k \times k$ matrices
$$\AM {\Phi^A} {M_k(C_0(\R^{2n}) \otimes A)} {M_k(\KK(L^2\R^n) \otimes A) \cong \KK(L^2(\R^n, A^k))}$$
by applying element-wise. Also extend $C_t^A : M_k(C_0(\R^n) \otimes A) \to \L(L^2(\R^n, A^k))$
in a similar manner.

The following is the most basic result which
relates the operator $D$, its principal symbol $\sigma(D)$, and this
asymptotic morphism.

\proclaim{Lemma 4.8} If $D$ is a formally self-adjoint elliptic differential $A$-operator of order one on $\R^n$ with  
constant coefficients, then $D$ is essential and for every $\phi \in C_0(\R^n)$ and $f \in C_0(\R)$, we have that 
$$\Lim \|\Phi^A_t(\phi f(\sigma(D))) - M_\phi f(t^{-1}\bar{D})\| = 0.$$
\endproclaim

\demo{Proof} The fact that $D$ is essentially self-adjoint follows from the Fourier identity
$$\widehat{D\eta}(\xi) = \sym(D)(\xi)\hat{\eta}(\xi)$$
since the total symbol $\sym(D)(\xi)$ is a self-adjoint matrix for $\xi \in T_x^*\R^n \cong \R^n$.
For any $\eta \in \Domain(D) = C_c^\infty(\R^n, A^k)$ we have the estimate
$$\|(1 + D^2)\eta\| \geq \|\eta\|$$
from which it easily follows that $(1 + D^2)$ has dense range. Thus, $D$ is essential.

We define $$[\phi f(\sigma(D))](x,\xi) = \phi(x) f(\sigma(D)(\xi))$$ via the $C^*$-functional
calculus.  Since $f \in C_0(\R)$ and $D$ is
elliptic, $f(\sigma(D)) \in M_k(C_0(\R^n) \otimes A)$ and so 
$$\phi f(\sigma(D)) = \phi \otimes f(\sigma(D)) \in M_k(C_0(\R^{2n}) \otimes A).$$

From the definition of $C^A_t : C_0(\R^n, A) \to \L(L^2(\R^n, A))$ we see that
$$\widehat{[C^A_t(f(\sigma(D)))\eta]}(\xi) = f(\sigma(D)(t^{-1}\xi))\hat{\eta}(\xi)$$
where $\eta \in L^2(\R^n,A^k)$. Let $D_1 = D - b$ be the first-order part of $D$. It follows by Lemma 4.7 that
$$\lim_{t \to \infty} \|\Phi^A_t(\phi f(\sigma(D))) - M_\phi f(t^{-1}\bar{D}_1)\| = 0$$ since $D_1$ has
constant coefficients and the total symbol of $t^{-1}D_1$ is
$$\sym(t^{-1}D_1)(\xi) = \sigma(D)(t^{-1}\xi) = t^{-1}\sigma(D)(\xi)$$
for all $t \geq 1$. Now $D - D_1 = b$ is order zero and so the result
follows by part 2 of Lemma 4.5. \qed \enddemo

\proclaim{Lemma 4.9} Let $K$ be a compact subset of  $\R^n$. Let $\phi \in
C_c^\infty(\R^n)$ with $\supp(\phi)$ contained in $K$ and let $f \in
C_0(\R)$. Let $D$ be as in the previous lemma. Then for every $\epsilon >
0$, there is a $\delta > 0$ such that if $$B = \sum_1^n c_j(x)\Dx j + d(x)$$
is an essential $A$-operator on $\R^n$ and the coefficients of
$D$ and $B$ differ in the uniform norm on $K$ by $\delta$, then 
$$\|M_{\phi}f(t^{-1}\bar{D}) - M_{\phi}f(t^{-1}\bar{B})\| < \epsilon$$
for $t$ large enough.
\endproclaim

\demo{Proof} The set $C$ of all $f \in C_0(\R)$ for which the Lemma holds is
clearly closed under addition (by the triangle inequality) and is closed
under operator adjoints. Suppose $\{f_n\}$ is a sequence in $C$ such that
$f_n \to f \in C_0(\R)$. Consider the following inequality $$\|M_\phi
f(t^{-1}\bar{D}) - M_\phi f(t^{-1}\bar{B})\| \leq 
 \ 2 \|\phi\| \|f - f_n\| + \|M_\phi f_n(t^{-1}\bar{D}) - M_\phi f_n(t^{-
1}\bar{B})\|$$ which follows by the triangle inequality and the Spectral
Theorem. Given $\epsilon > 0$, choose $n > 0$ such that $\|f - f_n\| <
\epsilon/(2 \|\phi\|)$. Choose $\delta > 0$ such that the Lemma is true for
$f_n$ with $\epsilon/2$. It follows that the Lemma is true for $f$ with
this $\delta$ and so $C$ is closed. Now suppose the Lemma is true for $f$
and $g$. Consider the inequality 
$$\aligned \|M_\phi(fg)(t^{-1}\bar{D}) &- M_\phi(fg)(t^{-1}\bar{B})\|
\\ & \leq 
\|M_\phi f(t^{-1}\bar{D}) g(t^{-1}\bar{D}) - M_\phi f(t^{-1}\bar{D})g(t^{-
1}\bar{B})\| \\
& \quad \quad  + \|M_\phi f(t^{-1}\bar{D}) g(t^{-1}\bar{B}) -
 M_\phi f(t^{-1}\bar{B})g(t^{-1}\bar{B})\| \\
& \leq 2 \|[M_\phi, f(t^{-1}\bar{D})]\| \|g\| + \|f\| \|M_\phi g(t^{-1}\bar{D})
- M_\phi g(t^{-1}\bar{B})\| \\
&\quad \quad + \|g\| \|M_\phi f(t^{-1}\bar{D}) - M_\phi f(t^{-1}\bar{B})\|.
\endaligned$$ It follows that the Lemma is true for $fg$ by part 1 of
Lemma 4.5. Thus, $C$ is a closed $C^*$-subalgebra of $C_0(\R)$.

By the Stone-Weierstrass Theorem, it suffices to prove that $C$ contains
the resolvent functions. Let $\epsilon > 0$ be given. Consider
$$\eqalignno{M_\phi r_\pm(t^{-1}\bar{D}) &- M_\phi r_\pm(t^{-1}\bar{B})
= M_{\phi}(r_\pm(t^{-1}\bar{D}) - r_\pm(t^{-1}\bar{B})) & \cr
&= t^{-1}M_{\phi}r_\pm(t^{-1}\bar{B})(B - D)r_\pm(t^{-1}\bar{D}) & \cr
&= t^{-1}[M_{\phi}, r_\pm(t^{-1}\bar{B})](B - D)r_\pm(t^{-1}\bar{D}) & (1)
\cr
& \quad \quad + t^{-1}r_\pm(t^{-1}\bar{B}) M_{\phi}(B - D)r_\pm(t^{-
1}\bar{D}) & \cr}$$
By Lemma 4.5 again, we have that
$\Lim \|[M_{\phi}, r_\pm(t^{-1}\bar{B})]\| = 0.$
From elliptic operator theory \cite{MF80}, $M_\phi (B -  D) : H^1(\R^n, A^k) \to L^2(K, A^k)$
is bounded by a constant times the sum of terms involving
$$\aligned & \sup\|\phi(x)(c_j(x) - a_j)\|,\ \sup\|\phi(x)(d(x) - b)\|, \\
 & \sup\|{{\partial}\over{\partial x_i}}{(\phi(x)(c_j(x)-a_j))}\|,\text{ and }
\sup\|{{\partial}\over{\partial x_i}} {(\phi(d(x)-b))}\|.\endaligned$$
The result now follows by choosing $t$ large enough so that the first term
in ($1$) above has norm less than $\epsilon/2$ and $\delta > 0$ small
enough so the second term has norm less than $\epsilon/2$. \qed \enddemo

Finally, we come to the main theorem of this sections which relates the
principal symbol $\sigma(D)$ and the $A$-index asymptotic morphism
$\{\Phi^A_t\}$ to the asymptotic morphism $\{\E^D_t\}$ determined by $D$
in Theorem 4.6.

\proclaim{Theorem 4.10} Let $D$ be an essential $A$-operator on $\R^n$. Then for all $f \in C_0(\R)$
and $\phi \in C_0(\R^n)$, we have $$\Lim \|\Phi^A_t(\phi f(\sigma(D))) - M_\phi
f(t^{-1}\bar{D})\| = 0.$$
\endproclaim

\demo{Proof} Consider the inequality
$$\aligned \|\Phi^A_t(\phi &f(\sigma(D))) - M_\phi f(t^{-1}\bar{D})\| \leq
\|\Phi^A_t(\phi f(\sigma(D))) - \Phi^A_t(\psi f(\sigma(D)))\| \\&+ \|\Phi^A_t(\psi
f(\sigma(D))) - M_\psi f(\sigma(D)))\| + \|M_\psi f(t^{-1}\bar{D}) - M_\phi
f(t^{-1}\bar{D})\|.
\endaligned$$ This implies the inequality
$$\aligned \limsup_t &\|\Phi^A_t(\phi f(\sigma(D))) - M_\phi f(t^{-1}\bar{D})\|
\leq \\&2\|\phi - \psi\| \|f\| + \limsup_t\|\Phi^A_t(\psi f(\sigma(D))) -
M_\psi f(t^{-1}\bar{D})\|.\endaligned$$
Hence, we may assume that $\phi$ is smooth and compactly supported.

Let $K = \supp(\phi) \subset \R^n$, which is compact, and let $\epsilon >
0$. Since $D$ is elliptic and  self-adjoint ($\sigma(D)(x,\xi)^* =
\sigma(D)(x,\xi)$), it follows that
$$\phi f(\sigma(D)) \in M_k(C_0(T^*\R^n) \otimes A).$$
For each $x \in K$, let $D^x$ be the constant coefficient operator defined
by ``freezing'' the coefficients of $D$ at $x$, i.e.,
$$D^x = \sum_1^n a_j^x \Dx{j} + b^x$$
where $a_j^x = a_j(x)$ and $b^x = b(x)$. Thus, the principal symbol of $D^x$
is $$\sigma(D^x)(\xi) = \sigma(D)(x,\xi).$$

By the compactness of $K$ and the finite covering dimension
\cite{Mun75} of $\R^n$, we may choose $\delta > 0$ small enough and
cover $K$ regularly by finitely many small balls 
$$B_\delta(x_1), \dots, B_\delta(x_N),\quad x_i \in K,$$
so that no more than $d = n + 1$  balls intersect at any point in $K$ and a
subordinate partition of unity $\{\phi_i\}$ such that
$$\aligned &\|M_{\phi_i} \Phi^A_t(\phi f(\sigma(D))) - M_{\phi_i} \Phi^A_t(\phi
f(\sigma(D^{x_i})))\| \leq \frac{\epsilon}{2d} \\
&\|M_{\phi_i} M_\phi f(\bar{D}^{x_i}) - M_{\phi_i}M_{\phi}f(t^{-
1}\bar{D})\| \leq \frac{\epsilon}{2d}\endaligned$$
for all $i$ and $t$ large enough (using the previous lemma).

We wish to show that $\limsup_t \|\Phi^A_t(\phi f(\sigma(D))) - M_\phi f(t^{-
1}\bar{D})\| \leq \epsilon$.

Partition $\{1,2,\dots,N\}$ into $d$ sets $C_k$ such that for all $i,j \in
C_k$ with $i \neq j$, we have
$$B_\delta(x_i) \cap B_\delta(x_j) = \emptyset.$$
Let $A_t = \Phi^A_t(\phi f(\sigma(D))) - M_\phi f(t^{-1}\bar{D})$. It follows
that
$$A_t  = \sum_1^N \phi_i A_t \sim_a \sum_1^N \phi_i^{1/2}A_t\phi_i^{1/2}
= \sum_{k=1}^{d} \sum_{i \in C_k} \phi_i^{1/2}A_t\phi_i^{1/2}$$

Hence, we see that:
$$\aligned 
& \limsup_t \|A_t\| \leq \limsup_t \sum_{k=1}^{d}\|\sum_{i \in
C_k} \phi_i^{1/2}A_t\phi_i^{1/2}\| \\
&= \limsup_t \sum_{k=1}^{d} \max_{i \in
C_k}\|\phi_i^{1/2}A_t\phi_i^{1/2}\| 
\leq d \limsup_t \max_i \|\phi_i^{1/2}A_t\phi_i^{1/2}\|. \endaligned$$
Since $\phi_i^{1/2}A_t\phi_i^{1/2} \sim_a \phi_i A_t$, we have
$\limsup_t \| \phi_i^{1/2}A_t\phi_i^{1/2}\| \leq \limsup_t
\|\phi_iA_t\|.$
But, for all $i$ and large enough values of $t$,
$$\aligned & \|\phi_iA_t\| 
= \|\phi_i \Phi^A_t(\phi f(\sigma(D))) - \phi_i\phi f(t^{-1}\bar{D})\| \\
&\leq \|\phi_i\Phi^A_t(\phi f(\sigma(D))) - \phi_i \Phi^A_t(\phi f(\sigma(D^{x_i})))\|
 + \|\phi_i\phi f(t^{-1}\bar{D}^{x_i}) - \phi_i \phi f(t^{-1}\bar{D})\| \\
&\leq \frac{\epsilon}{2d} + \frac{\epsilon}{2d} < \epsilon \endaligned$$
since $\Phi^A_t(\phi f(\sigma(D^{x_i}))) \sim_a M_\phi f(t^{-1}\bar{D}^{x_i})$ by Lemma 4.8. The
result follows. \qed \enddemo

\head {\bf 5. The Exact Mishchenko-Fomenko Index Theorem} \endhead

Let $A$ be a unital $C^*$-algebra. Let $M$ be a smooth closed Riemannian manifold and let $D : C^\infty(E) \to
C^\infty(F)$ be an elliptic differential $A$-operator of order one on $M$. The following is the main example in the theory,
which we will need in the proof of the index theorem. For the basic theory of (pseudo)differential $A$-operators
on vector $A$-bundles see \cite{MF80,Mis78}.

\demo{\bf Example 5.1} Let $S$ be a smooth Clifford bundle on $M$ with
Clifford multiplication $c : T^*M \to \End(S)$. Let $D : C^\infty(S)
\to C^\infty(S)$ be the associated Dirac-type operator \cite{LM89} defined via a compatible
connection $\nabla^S: C^\infty(S) \to C^\infty(T^*M \otimes S)$. 
Let $E$ be a smooth vector $A$-bundle on $M$ also equipped with
a connection $\nabla^E : C^\infty(E) \to C^\infty(T^*M \otimes E)$. Let $\nabla =
\nabla^S \otimes 1 + 1 \otimes \nabla^E$ denote the tensor product
connection on the vector $A$-bundle $S \otimes E$. The associated
generalized Dirac operator $D_E : C^\infty(S \otimes E) \to C^\infty(S
\otimes E)$ is an elliptic differential $A$-operator of order one on $M$.
In local coordinates,
$$D_E = \sum_{j=1}^n c(e_j) \nabla_{e_j}$$
where $\{e_j\}$ is a local orthonormal frame for $TM \cong T^*M$.
\enddemo 

Since $D$ is elliptic, we may assume $F = E$. Put a smooth Hermitian $A$-metric on the vector $A$-bundle $E$. Consider $D$ as a densely-defined,
unbounded $A$-operator on the Hilbert $A$-module $L^2(E)$ with dense $\Domain(D) = C^\infty(E)$. If $D$ is
formally self-adjoint, then the closure of the graph of $D$ in $L^2(E)
\oplus L^2(E)$ is the graph of an unbounded $A$-operator $\bar{D}$, called
the {\it closure} of $D$.  The basic elliptic estimate \cite{MF80} shows that
$\Domain(\bar{D})$ is the Sobolev space $H^1(E)$.
We refer the reader to Definition 4.2 for the following.

\proclaim{Lemma 5.2} If $D$ is a formally self-adjoint elliptic differential $A$-operator of order one on $M$, then $D$ is essential and $\bar{D}$ has compact
resolvents. Moreover, the restriction $D|_U : C^\infty(U,E) \to C^\infty(U,E)$ to any open subset $U$ of $M$ is also essential. 
\endproclaim

\demo{Proof} The fact that $\bar{D}$ is self-adjoint follows from ellipticity and the existence of a parametrix \cite{MF80} 
as in the classical case. The estimate \cite{MP93} $$\|(1+ D^2)\eta\| \geq \|\eta\|$$ shows that $(1 + D^2)$ has dense range.
From the identity $(D + i)(D - i) = (1 + D^2)$ we have that $(D \pm i)$ have dense range and so $D$ is essential.
The generalized Rellich Lemma and basic elliptic estimate then prove that $(\bar{D} \pm i)^{-1}$ are compact $A$-operators
on $L^2(E)$. Note that each resolvent is the adjoint of the other. Using the fact that $\supp(D \eta) \subset \supp(\eta)$ for {\it differential}
operators and the decomposition $L^2(E) \cong L^2(U, E) \oplus L^2(M \backslash U, E)$, the last statement follows. \hfill \qed
\enddemo 

Consider the formally self-adjoint $A$-operator 
$$\b{D} = \pmatrix 0 & D^t \\ D & 0 \endpmatrix : C^\infty(E \oplus E) \to
C^\infty(E \oplus E)$$ where $D^t$ is the formal adjoint of $D$.  
The principal symbol of $D$ is the self-adjoint
homomorphism
$$\bsigma = \sigma(\b D) = \pmatrix 0 & \sigma(D)^* \\ \sigma(D) & 0 \endpmatrix : \pi^*(E
\oplus E) \to \pi^*(E \oplus E),$$ where $\sigma(D) : \pi^*E \to \pi^*E$ is the
principal symbol of $D$.

\proclaim{Lemma 5.3} The resolvents $$(\bsigma \pm i)^{-1} : \pi^*(E \oplus E) \to
\pi^*(E \oplus E)$$ are $A$-homomorphisms which vanish at infinity on
$T^*M$ in the operator norm induced by the Hermitian $A$-metrics on $E$.
\endproclaim

\demo{Proof} Follows from homogeneity $\bsigma(x,t\xi) = t \bsigma(x,\xi)$ and
ellipticity. \qed \enddemo

Form the Cayley transform ~\cite{Qui88}
$$\bold u = (\bsigma + i)(\bsigma - i)^{-1} = 1 + 2i(\bsigma - i)^{-1}.$$ By
complementing the vector $A$-bundle $E$, we may embed
$\pi^*(E \oplus E)$ in a trivial $A$-bundle $$\b A = T^*M \times (A^{n}
\oplus A^{n}).$$ 
Now extend the automorphism $\bold u$ to the $A$-bundle $\b A$ by
defining it to be equal to the identity on the complement of $\pi^* E \oplus
\pi^* E$ in $\b A$. From the lemma above, it follows that $\bold u$
extends continuously to the trivial $A$-bundle on the one-point
compactification $(T^*M)^+$ by setting $\bold u(\infty) = I$.

Let $$\bold \epsilon = \pmatrix 1 & 0 \\ 0 & -1 \endpmatrix$$ be the
grading of the trivial $A$-bundle $(T^*M)^+ \times (A^{n} \oplus A^{n})$.
Since $\bold \epsilon \bsigma = - \bsigma \bold \epsilon$ it follows that 
$(\bold u \bold \epsilon)^2 = 1.$ A simple calculation also shows that
$(\bold u \bold \epsilon)^* = \bold u \bold \epsilon$ is self adjoint. We
also have obviously that $\epsilon^* = \epsilon$ and $\epsilon^2 = 1$. 

Recall that for every self adjoint involution $x$ there is an associated
projection $p(x) = \frac12(x + 1).$ In our case, we obtain two
projection-valued functions $p(\bold \epsilon)$ and $p(\bold u \bold
\epsilon)$ on $(T^*M)^+$ which are equal at infinity. Both define elements
in $$K_0(C(T^*M^+) \otimes A) = K^0_A(T^*M^+)$$ and so their difference
defines an element in
$$\sigma_D = [p(\bold \epsilon)] - [p(\bold \epsilon \bold u)] \in
K_0(C_0(T^*M) \otimes A) = K^0_A(T^*M).$$ This is the {\it symbol class}
of the elliptic $A$-operator $D$ as constructed in \cite{Hig93} for $A = \C$.
(See also Quillen \cite{Qui88}.)

\proclaim{Lemma 5.4} $\sigma_D = [\sigma(D)] \in K^0_A(T^*M).$ \endproclaim

For each $t \geq 1$, we can form the compact $A$-operators \footnote{We will write $\b D$ instead of the closure $\bar{\b D}$ from now on.}
$$f(t^{-1}\b D) : L^2(M, E \oplus E) \to L^2(M, E \oplus E)$$
by Lemma 5.3 and the functional calculus.  We can extend $f(t^{-1}\b D)$ to an operator on the Hilbert $A$-module
$L^2(M, A^{n} \oplus A^{n})$  by defining it to be zero on the complement of
the $A$-submodule $L^2(M,E \oplus E)$ (which exists by complementing
$E$ as above) in $L^2(M, A^n \oplus A^n) = L^2(M, A^{2n})$. 
This then defines a continuous family of $*$-homomorphisms
$$C_0(\R) \to \KK(L^2(M, A^{2n})) : f \mapsto f(t^{-1}\b D).$$

Likewise, we may apply $f$ to the symbol $\bsigma$ of $\b D$. This $A$-homomorphism 
$f(\bsigma)$ of the vector $A$-bundle $\pi^*(E \oplus F)$ also vanishes at infinity (compare Lemma 5.3). 
We may then extend $f(\bsigma)$ to an  $A$-homomorphism of the trivial $A$-bundle $\b A$ by setting it equal to zero
on the complement of $\pi^*(E \oplus F)$.

Let $\amn {\Phi^{M,A}} {M_{2n}(C_0(T^*M) \otimes A)} {\KK(L^2(M,A^{2n}))}$ denote the $A$-index asymptotic morphism for $M$
from Appendix $A$. By thinking of $f(\bsigma)$ as a matrix of $A$-valued functions on $T^*M$ vanishing at infinity,
we obtain a $*$-homomorphism
$$C_0(\R) \to C_0(T^*M, M_{2n}(A)) = M_{2n}(C_0(T^*M) \otimes A) : f \mapsto f(\bsigma).$$
The following result relates the spectral theory of $D$ to the principal symbol $\sigma(\b D)$ and is
central to the index theorem.

\proclaim{Lemma 5.5} If $D$ is an elliptic differential $A$-operator of order one on $M$ with symbol
$\sigma$, then for every $f \in C_0(\R)$,
$$\Lim \|\Phi_t^{M,A}(f(\sigma (\b D))) - f(t^{-1}\b D)\| = 0$$ \endproclaim

\demo{Proof} By complementing the vector $A$-bundle $E$ we may assume that $E = M \times A^n$ is actually trivial.
Let $\{U_j\}_1^m$ be an open cover of $M$ by coordinate charts. Let $\b D_j = \b D|_{U_j}$ be the restriction 
to $U_j$. Choose a smooth partition of unity $\{\rho_j^2\}$ subordinate to the cover $\{U_j\}$. Let 
$$\{\Phi_t^j\} = \{\Phi_t^{M,A}|_{U_j}\} : M_{2n}(C_0(T^*U_j) \otimes A) \to \KK(L^2(U_j,A^{2n}))$$
denote the restriction to $U_j$ as in Corollary A.11. By Lemma 5.2 and Theorem 4.10, we have for each $1 \leq j \leq m$
$$\lim_{t \to \infty} \| \Phi_t^j(\rho_j f(\sigma(\b D_j))\rho_j)  -  M_{\rho_j}f(t^{-1}\b D_j) M_{\rho_j}\| = 0.$$
Consider the following asymptotic equivalences:
$$\aligned \Phi_t^{M,A}(f(\sigma(\b D))) &= \Phi_t^{M,A}(\sum_1^m \rho_j f(\sigma(\b D)) \rho_j) \sim_a \sum_1^m \Phi_t^{M,A}(\rho_j f(\sigma(\b D)) \rho_j)\\
& =  \sum_1^n \Phi_t^{M,A}(\rho_j f(\sigma(\b D_j)) \rho_j) \sim_a \sum_1^n \Phi_t^j(\rho_j f(\sigma(\b D_j)) \rho_j)\\
&\sim_a \sum_1^m M_{\rho_j} f(t^{-1}\b D_j) M_{\rho_j} =  \sum_1^m M_{\rho_j} f(t^{-1}\b D) M_{\rho_j} \\
&\sim_a \sum_1^m  M_{\rho_j} M_{\rho_j} f(t^{-1} \b D) = f(t^{-1} \b D) \endaligned$$
This completes the proof. \qed \enddemo

\proclaim{Proposition 5.6} Let $D$ be an elliptic differential $A$-operator of order one on the
smooth closed manifold $M$. The analytic and topological indices of $D$
are equal, that is, $$\Index_a(D) = \Index_t(D) \in K_0(A).$$
\endproclaim

\demo{Proof} For all $t > 0$ form the Cayley transform
$$U_t = (t^{-1}\b D + i)(t^{-1}\b D - i)^{-1} = I + 2i(t^{-1}\b D - i)^{-1}$$
which defines a unitary $A$-operator on $L^2(M, A^{2n})$ which is 
equal to the identity on the complement of $L^2(M, E \oplus E)$.
It follows from Lemma  5.5 that
$$\Lim \|U_t - \Phi^{M,A}_t(\bold u) \| = 0$$ (where we extended the result by
adjoining a unit.) It follows from this that $\Phi^{M,A}_*([p(\epsilon \bold u)]) = [p(\epsilon U_1)]$ and so the morphism index 
can be written as
$$\Ind_m(M,\sigma(D)) = \Phi^{M,A}_*([\sigma(D)]) = \Phi^{M,A}_*(\sigma_D) = [p(\bold
\epsilon)] - [p(\bold \epsilon U_1)] \in K_0(A).$$
Let us now ponder what happens as $t \to 0$. By homotopy invariance, we
may assume that the kernel and cokernel of $D$ are finitely generated and
projective $A$-modules and that $\b D$ is a Fredholm $A$-operator with
closed range. Thus, we have the Hilbert $A$-module decomposition
$$L^2(M, A^n \oplus A^n) = \Ker(\b D) \oplus \H  = \Ker(D) \oplus \Ker(D^*)
\oplus \H.$$
For all $t$, we have that
$$U_t|_{\Ker(\b D)} = I + 2i(0 - iI)^{-1} = I - 2I = -I$$ on the kernel of $\b
D$. Thus,
$$\bold \epsilon U_t|_{\Ker(\b D)} = \pmatrix -1 & 0 \\ 0 & 1
\endpmatrix$$ with respect to the decomposition $\Ker(\b D) = \Ker(D)
\oplus \Ker(D^*)$.
On the orthogonal complement of $\Ker(\b D)$, $U_t$ converges in norm to
the identity as $t \to \infty$, since $\b D$ is invertible off $\Ker(\b D)$
and has a gap in the spectrum. Thus, $U_t$ converges in norm to $U_0$
which is $-I$ on $\Ker(\b D)$ and $I$ on the complement. A simple
calculation then shows that
$$[p(\bold \epsilon)] - [p(\bold \epsilon U_1)] = [p(\bold \epsilon)] - [p(\bold \epsilon U_0)] = [P_{\Ker(D)}] -
[P_{\Ker(D^*)}]$$ where $P_{\Ker(D)}$ and $P_{\Ker(D^*)}$ are the
projections onto the kernel and cokernel of $D$. Therefore, we have by Theorem 3.13 that
$$\Index_t(D) = \Ind_t(M,\sigma(D)) = \Ind_m(M,\sigma(D)) = \Index_a(D).  \qed $$ 
\enddemo

This establishes the exact version of the Mishchenko-Fomenko Index Theorem for first-order elliptic
differential $A$-operators on {\it arbitrary} smooth closed manifolds. It is not known if these techniques
can be extended verbatim to higher-order differential (or even pseudodifferential)
$A$-operators since the results in Section 4, especially the estimates in Lemma 4.9, were specific
to first-order operators via the symbol identity. However, using Kasparov's $KK$-theory \cite{Kas81},
we can derive the full Mishchenko-Fomenko Index Theorem for \spinc manifolds using twisted
Dirac operators \cite{LM89}. See the helpful discussions in Section 24 of Blackadar \cite{Bla86}.

Let $M$ be a smooth compact \spinc-manifold without
boundary and let $D$ be the Dirac operator of $M$. Let $\dolbeaut$ denote
the Dolbeaut operator on the cotangent bundle $T^*M$. These two elliptic
operators determine $KK$-classes
$$\aligned [D] &\in KK^0(C(M),\C) \\
	    [\dolbeaut] &\in KK^0(C_0(T^*M),\C)
\endaligned$$
The Thom isomorphism determines a class
$$x = [\Psi^{T^*M}] \in KK^0(C(M),C_0(T^*M)).$$ 
This element relates these  two elliptic classes (Lemma 24.5.1 \cite{Bla86}) by 
$$x \otimes_{C_0(T^*M)} [\dolbeaut] = [D],$$ where $\otimes$ denotes
Kasparov's intersection product (and corresponds to composition of
asymptotic morphisms).

Let $P : C^\infty(E_1) \to C^\infty(E_2)$ be an elliptic pseudodifferential $A$-operator of order $m$ on $M$.
$P$ determines a $KK$-class $[P] \in KK^0(C(M),A)$. (Compare Theorem 4.6). The principal symbol $\sigma(P)$ of $P$ determines two $KK$-classes
$$\aligned &[[\sigma(P)]] \in KK^0(C(M),C_0(T^*M) \otimes A) \\
	   &[\sigma(P)] \in KK^0(\C,C_0(T^*M) \otimes A) = K^0_A(T^*M).
\endaligned$$ They are related by $f^M_*[[\sigma(P)]] = [\sigma(P)]$ where $f^M : M \to \{pt\}$ is the collapsing map.

\proclaim{Lemma 5.7} {\rm (Theorem 5 \cite{Kas84}) $[P] = [[\sigma(P)]] \otimes_{C_0(T^*M)} [\dolbeaut].$
\endproclaim

The {\it analytic index} of $D$ is given by the formula
$$\Index_a(D) = f^M_*([P]) \in KK^0(\C,A) = K_0(A).$$
The {\it topological index} of $P$ is defined to be $$\Index_t(P) = \Ind_t(M, \sigma(P)) \in K_0(A)$$
where $\sigma(P) : \pi^*(E_1) \to \pi^*(E_2)$ is the principal symbol of $P$.

\proclaim{Theorem 5.8 (Exact Mishchenko-Fomenko Index Theorem)} If $P$ is an elliptic pseudodifferential $A$-operator
on the smooth closed $\text{spin}^c$-manifold $M$, then the analytic and
topological indices of $P$ are equal, that is,
$$\Index_a(P) = \Index_t(P) \in K_0(A).$$
\endproclaim

\demo{Proof} We may suppose $M$ is even-dimensional.
(If $\dim(M) = n$ is odd, then replace $M$ by $M \times \Bbb T^N$.) 

By the Thom isomorphism, there is a vector $A$-bundle $E$ on $M$ such that 
$$[E] \otimes_{C(M)} x = [\sigma(P)],$$ 
where $[E] \in K^0_A(M) = KK^0(\C,C(M) \otimes A)$ is the associated $K$-theory class. It also
determines a class $[[E]] \in KK^0(C(M), C(M) \otimes A)$, where
$f^M_*[[E]] = [E]$. 

Let $D_E$ denote the Dirac operator of $M$ twisted by the vector $A$-bundle $E$
(as in Example 5.1).
Note that $D_E$ is a first-order elliptic differential $A$-operator on $M$
and so $$\Index_a(D_E) = \Index_t(D_E)$$ by the previous theorem.  By Lemma 24.5.3 \cite{Bla86},
$$[D_E] = [[E]] \otimes_{C(M)} [D].$$

The $A$-index asymptotic morphism of $M$ determines a class
$$[\Phi_t^{M,A}]  \in E^0(C_0(T^*M) \otimes A, A) \cong KK^0(C_0(T^*M) \otimes A, A).$$

We then compute that:
$$\aligned 
\Index_a(P) & = f^M_*[P] = f^M_*([[\sigma(P)]] \otimes_{C_0(T^*M)} [\dolbeaut]) \\
& = [\sigma(P)] \otimes_{C_0(T^*M)} [\dolbeaut] = [E] \otimes_{C(M)} x \otimes_{C_0(T^*M)} [\dolbeaut] \\
& = [E] \otimes_{C(M)} [D] = f^M_*([[E]] \otimes_{C(M)} [D]) = f^M_*[D_E] \\
& = \Index_a(D_E) = \Index_t(D_E) = \Ind_m(M, \sigma(D_E))\\
& = [\sigma(D_E)] \otimes_{C_0(T^*M) \otimes A} [\Phi_t^{M,A}] = [E] \otimes_{C(M)} x \otimes_{C_0(T^*M) \otimes A}
[\Phi_t^{M,A}] \\
& = [\sigma(P)] \otimes_{C_0(T^*M) \otimes A} [\Phi_t^{M,A}] = \Ind_m(M, \sigma(P)) = \Index_t(P)
\endaligned$$
and we are done. \qed \enddemo


\head {\bf Appendix A: The Index Asymptotic Morphism} \endhead 

We will construct for each smooth Riemannian manifold $M$ without
boundary and $C^*$-algebra $A$ a natural asymptotic morphism \cite{CH89}
$$\AM {\Phi^A} {C_0(T^*M) \otimes A} {\KK(L^2M) \otimes A}.$$
When $A = \C$ is the complex numbers, this asymptotic morphism is equivalent to the one
originally described by Connes and Higson \cite{CH89} and elaborated on by Higson \cite{Hig93}, although we
are presenting it from a different viewpoint.

Fix an integer $n \geq 0$. Denote by $\L(L^2(\R^n))$ the $C^*$-algebra of bounded linear operators on
$L^2\R^n$ and let $C_0(\R^n)$ denote the $C^*$-algebra of continuous functions vanishing at infinity.
Let $C^*(\R^n)$ be the $C^*$-algebra completion of the convolution algebra $L^1(\R^n)$ in the operator norm of
$\L(L^2\R^n)$. The Fourier Transform gives a continuous algebra
homomorphism $L^1(\R^n) @>\wedge>> C_0(\R^n)$ which extends to a
$C^*$-algebra isomorphism $C^*(\R^n) @>\wedge>> C_0(\R^n)$ (called the
Gel'fand Transform \cite{Rud91}.) For $t \geq 1$ and $f \in C_0(\R^n)$
with $\check{f} \in L^1(\R^n)$, we define $C_t(f) : L^2(\R^n) \to L^2(\R^n)$ by
$$[C_t(f)\phi](x) = \check{f_t}*\phi(x)$$
where $\phi \in L^2(\R^n)$ and $f_t(x) =_{\text{def}} f(t^{-1}x)$ for $x \in \R^n$.

\proclaim{Lemma A.1} $C_t$ extends to a continuous family of $*$-monomorphisms
$$\AM {C} {C_0(\R^n)} {\L(L^2\R^n)}.$$
\endproclaim

For $f \in C_0(\R^n)$ let $M_f \in \Cal B(L^2\R^n)$ denote multiplication by $f$. This
defines a $*$-homomorphism $M : C_0(\R^n) \to \Cal B(L^2\R^n)$.

\proclaim{Proposition A.2} For all functions $f, g \in C_0(\R^n)$,
$$\Lim \|M_f C_t(g) - C_t(g) M_f \| = 0.$$
\endproclaim

\demo{Proof} We may assume that $f$ is
real-valued and compactly supported. The set
$\Cal A$ of all $g$ which satisfy the conclusion above forms a $C^*$-subalgebra of $C_0(\R)$.

\noindent (a.) If $n = 1$, then we need only to check that $\Cal A$ contains
the generators $r_\pm(x) = (x \pm i)^{-1}$ of $C_0(\R)$. We have that
$$C_t(r_\pm) = r_\pm(t^{-1} D) = (t^{-1}D \pm i)^{-1}$$
where $D = -\sqrt{-1} \frac{d}{dx}$ is the (closure of the) Dirac operator.
The result now follows from the commutator identity
$$[M_f, r_\pm(t^{-1}D)] = r_\pm(t^{-1}D)[t^{-1}D, M_f]r_\pm(t^{-1}D)$$
and noting that the bounded operators $[t^{-1}D,M_f] = t^{-1}
[D,M_f] \to 0$ as $t \to \infty$.

\noindent (b.) If $n > 1$, then $C_0(\R^n) \cong C_0(\R) \otimes C_0(\R)
\otimes \cdots \otimes C_0(\R)$. If $g \in C_0(\R^n)$ has the form 
$g(\xi) = g_1(\xi_1) \otimes \cdots g_n(\xi_n)$
then it follows that $C_t(g) = g_1(t^{-1}D_1) \cdots g_n(t^{-1}D_n)$
where $D_j = \D j$. Apply part (a.) inductively. \qed \enddemo

The previous result says that the continuous family $\{C_t\}$ and the constant family $\{M_t = M\}$
``asymptotically commute'' with each other.

\demo{\bf Definition A.3} Using the canonical identification
$T^*\R^n \cong \R^n \times \R^n$, which identifies $C_0(T^*\R^n) \cong C_0(\R^n) \otimes C_0(\R^n)$,
we define the asymptotic morphism
$$\AM {\Phi} {C_0(T^*\R^n)} {\L(L^2\R^n)}$$
to be (induced by) the tensor product $f \otimes g \mapsto M_f \circ C_t(g)$. This asymptotic morphism is only
well-defined up to asymptotic equivalence. \enddemo


We now want to show that the image of each $\Phi_t$ actually lies in the $C^*$-algebra
of compact operators $\KK(L^2\R^n)$.

\proclaim{Theorem A.4} $\AM {\Phi} {C_0(T^*\R^n)} {\KK(L^2\R^n)}$
\endproclaim

\demo{Proof}  Let $F = f \otimes g$ where $f, g$ are {\it compactly supported}.
For $\eta \in L^2(\R^n)$:
$$[\Phi_t(F)\eta](x) = [M_f\circ C_t(g)\eta](x)= \int_{\R^n}
k^F_t(x,y)\eta(y)\ dy,$$
where the kernel
$$k^F_t(x,y) = \frac{1}{(2\pi)^{n/2}}f(x)\check{g_t}(x-y) = 
(\frac{t}{2\pi})^n \int_{\R^n} f(x)g(\xi)e^{\sqrt{-1}t(x-y)\xi}\ d\xi$$
is in $L^2(\R^n \times \R^n)$. Thus, $\Phi_t(F)$ is an integral operator
with square-integrable kernel and so is a compact operator on $L^2(\R^n)$. Now take a norm limit. \qed
\enddemo

\proclaim{Lemma A.5} For all $\rho\in C_0(\R^n)$ and $F \in C_0(T^*\R^n)$, we
have that:
\roster
\item $\Phi_t(\rho F) = M_\rho \Phi_t(F)$, for all $t \geq 1$;
\item $\Lim \|M_\rho \Phi_t(f) - \Phi_t(f)M_\rho\| = 0.$
\endroster
\endproclaim

\demo{Proof} Part (1) follows from the fact that $M_{\rho f} = M_\rho M_f$. The
second follows from a similar argument in Proposition A.2.
\qed \enddemo

This asymptotic morphism also ``restricts'' to open subsets. Let $U$ be an
open subset of $\R^n$. Then $T^*U \cong U \times \R^n \subset T^*\R^n$ is
also open and we have the natural inclusion $C_0(T^*U) \subset
C_0(T^*\R^n)$ and decomposition $L^2(\R^n) = L^2(U) \oplus L^2(\R^n
\setminus U)$. 

\proclaim{Lemma A.6} For all $t \geq 1$, $\Phi_t|_{C_0(T^*U)} : C_0(T^*U) \to \KK(L^2U)$.
\endproclaim

\demo{Proof} If $F = f \otimes g \in C_0(T^*U) = C_0(U) \otimes C_0(\R^n)$, then
for all $\eta \in L^2(U)$, it follows that
$$\supp (\Phi_t(F)\eta) = \supp (M_f C_t(g) \eta) \subset \supp (f)
\subset U.$$
The result now follows by an approximation argument since each $\Phi_t$ is linear.
\qed \enddemo

Suppose that $\psi : U \to W$ is a diffeomorphism of open subsets of
$\R^n$.  Denote by $\hat{\psi} : T^*U \to T^*W$ the induced diffeomorphism
of cotangent bundles which is defined by
$$\hat{\psi}(x,\xi) = (\psi(x), d\psi^{-t}\xi),$$
where $d\psi : TU \to TW$ denotes the derivative of $\psi$, mapping
tangent vectors at $x$ to tangent vectors at $\psi(x)$, and $d\psi^{-t}$
denotes the inverse of the transpose, mapping cotangent vectors at $x$ to
cotangent vectors at $\psi(x)$. Let $T_{\psi} : L^2(W) \to L^2(U)$ denote
the induced unitary isomorphism of Hilbert spaces defined by
$$T_{\psi}\eta(x) = \eta(\psi(x))\  J^{1/2}(x),\ \eta \in L^2(W),$$
where $J(x)$ denotes the absolute value of the Jacobian determinant of $\psi$ at $x \in U$.

\proclaim{Proposition A.7}{\rm (Lemma 8.7 \cite{Hig93})} If a function $F \in C_0(T^*\R^n)$ has support in the
open subset $T^*W$, then
$$\Lim \|\Phi_t(F \circ \hat{\psi}) - T_{\psi}\Phi_t(F)T_{\psi}^{-1}\| =
0.$$
\endproclaim

Now let $M$ be a smooth Riemannian $n$-manifold without boundary.
Cover $M$ with open charts $U_\alpha$, each diffeomorphic to an open
subset $W_\alpha$ of $\R^n$ via the diffeomorphism $\psi_\alpha :
U_\alpha \to W_\alpha$. Let $\check{\psi}_\alpha = ({\psi_\alpha}_*^t) :
T^*W_\alpha \to T^*U_\alpha$ be the induced diffeomorphism of cotangent
bundles and let $T_\alpha : L^2(W_\alpha) \to L^2(U_\alpha)$ be the
induced unitary isomorphism of Hilbert spaces. Define $\Phi^\alpha_t :
C_0(T^*U_\alpha) \to \KK(L^2U_\alpha)$ by the following
$$\Phi^\alpha_t(f) = T_\alpha \Phi_t(f \circ \check{\psi}_\alpha)
T_\alpha^{-1}.$$ (We include the appropriate Radon-Nikodym derivative in
the definition of $T_\alpha$.) That is, we define $\Phi^\alpha_t$ so that
the following diagram commutes asymptotically:
$$\CD 
C_0(T^*U_\alpha) @>\Phi^\alpha_t>> \KK(L^2U_\alpha) \\
@VV{\cong}V                                  @AA{\cong}A  \\
C_0(T^*W_\alpha) @>{\Phi_t|_{W_\alpha}}>> \KK(L^2W_\alpha) \\
@VVV                                  @VVV  \\
C_0(T^*\R^n)     @>\Phi_t>>        \KK(L^2\R^n) 
\endCD$$
where the bottom diagram commutes by Lemma A.6.

\proclaim{Corollary A.8} If $f \in C_0(T^*M)$ has support in $U_\alpha \cap U_\beta$, then
$$\Lim \|\Phi^\alpha_t(f) - \Phi^\beta_t(f)\| = 0$$
in the operator norm on $L^2(M)$.
\endproclaim

Now let $\{\rho^2_\alpha\}$ be a smooth partition of unity subordinate to
the open cover $\{U_\alpha\}$. Define the family of functions
$$\AM {\Phi^M} {C_0(T^*M)} {\KK(L^2M)}$$
in the following way:
$$\Phi_t^M(f)\eta = \sum_{\alpha} \Phi^{\alpha}_t(\rho_{\alpha} f) (\rho_{\alpha} \eta),$$
where $f \in C_0(T^*M)$ and $\eta \in L^2(M)$.

\proclaim{Theorem A.9} $\AM {\Phi^M} {C_0(T^*M)} {\KK(L^2M)}$ is an asymptotic
morphism which is asymptotically independent of the choice of open cover
and partition of unity. Moreover, if $\psi : U \to W$ is a diffeomorphism
from an open subset in $M$ to an open subset in $\R^n$ then for all $f \in C_c^{\infty}(T^*W)$
$$\Lim \|\Phi_t^M(f \circ \check{\psi}) - T_{\psi}\Phi_t(f)T^{-1}_{\psi}\| = 0$$
\endproclaim

\demo{\bf Definition A.10} Let $A$ be a $C^*$-algebra. Define the 
{\it $A$-index asymptotic morphism for $M$},
$$\AM {\Phi^{M,A}} {C_0(T^*M) \otimes A} {\KK(L^2M) \otimes A},$$
to be the (asymptotic morphism) tensor product
of the index asymptotic morphism $\amn {\Phi^M} {C_0(T^*M)} {\KK(L^2M)}$ 
above with the identity morphism $id_A : A \to A$.
\enddemo

The following properties follow easily from the previous results. 

\proclaim{Corollary A.11} Let $U \subset M$ be an open subset of $M$. Then
$C_0(T^*U)$ is an ideal in $C_0(T^*M)$ and $\KK(L^2U) \hookrightarrow
\KK(L^2M)$. The following diagram commutes asymptotically
$$\CD
C_0(T^*M) \otimes A @>{\{\Phi^{M,A}_t\}}>> \KK(L^2M) \otimes A \\
@AAA                                     @AAA \\
C_0(T^*U) \otimes A @>{\{\Phi^{M,A}_t|_U\}}>> \KK(L^2U) \otimes A 
\endCD$$
\endproclaim

\proclaim{Corollary A.12} Let $\psi : M \to M$ be a smooth diffeomorphism. The
following diagram commutes asymptotically:
$$\CD
C_0(T^*M) \otimes A @>{\{\Phi^{M,A}_t\}}>> \KK(L^2M) \otimes A \\
@AA{\psi^*\otimes id_A}A           @AA{Ad(\psi^*)\otimes id_A}A \\
C_0(T^*M) \otimes A @>{\{\Phi^{M,A}_t\}}>> \KK(L^2M) \otimes A 
\endCD$$
\endproclaim


\head {\bf Appendix B: Bott Periodicity and Thom Isomorphism} \endhead

Let $A$ be a $C^*$-algebra. In this section, 
we prove that the asymptotic morphism
$$\AM {\Phi} {C_0(\R^{2n}) \otimes A} {\KK(L^2\R^n) \otimes A},$$
constructed in Appendix A, induces the (inverse of the) Bott Periodicity 
isomorphism $K_0(C_0(\R^{2n}) \otimes A) \cong K_0(A)$. The proof is a 
direct generalization of Atiyah's functorial proof of Bott Periodicity in 
topological $K$-theory \cite{Ati68}.

Suppose $A$ has a unit. If $B$ is another unital $C^*$-algebra, there is a
well-defined map
$$\aligned \mu : K_0(A) \otimes K_0(B) &\to K_0(A \otimes B) \\
[p] \otimes [q] &\mapsto [p \otimes q] \endaligned$$ where $p$ is a
projection over $A$ and $q$ is a projection over $B$. (If $A$ or $B$ has no
unit, then adjoin one and note that the above map for the unitized algebras
restricts as needed.)

\proclaim{Theorem B.1} Suppose for each $C^*$-algebra $A$ there is a
homomorphism $$\alpha_A : K_0(C_0(T^*\R^n) \otimes A) \to K_0(A)$$
which satisfies the following properties:
\roster
\item $\alpha_A$ is natural in $A$, i.e. for every morphism $\psi : A \to B$
the following diagram commutes:
$$\CD
{K_0(C_0(T^*\R^n)\otimes A)} @>{\alpha_A}>>   {K_0(A)} \\
   @V{\psi_*}VV                           @V{\psi_*}VV \\
{K_0(C_0(T^*\R^n)\otimes B)} @>{\alpha_B}>>   {K_0(B)} \\ 
\endCD$$
\item For any $C^*$-algebra $B$, the diagram below commutes:
$$\CD
{K_0(C_0(T^*\R^n) \otimes A) \otimes K_0(B)} @>{\alpha_A \otimes
id_B}>> {K_0(A) \otimes K_0(B)} \\
@V{\mu}VV            @V{\mu}VV \\
{K_0(C_0(T^*\R^n)\otimes A \otimes B)} @>{\alpha_{A\otimes B}}>>
{K_0(A \otimes B)}
\endCD$$
\item There is an element $b \in K_0(C_0(T^*\R^n))$ such that $\alpha_{\C}(b) = 1
\in K_0(\C) = \Z$.
\endroster
Then $\alpha_A$ is an isomorphism for all $A$.
\endproclaim

\demo{Proof} The inverse map $\beta_A : K_0(A) \to K_0(C_0(T^*\R^n) \otimes A)$ is defined by the formula
$$\beta_A(x) = \mu (b \otimes x),$$
for $x \in K_0(A).$ Now use properties (1.) - (3.) to verify this. \hfill \qed \enddemo

\proclaim{Lemma B.2} $\Phi^A_*$ satisfies properties {\rm (1)} and {\rm (2)}. \endproclaim

\demo{Proof} Property (1) follows from the fact, since $\Phi_t$ is linear, the
diagram
$$\CD
{C_0(T^*\R^n) \odot A} @>{\Phi_t \otimes id_A}>> \KK \odot A \\
@V{id \otimes \psi}VV                           @V{id \otimes \psi}VV \\
{C_0(T^*\R^n) \odot B} @>{\Phi_t \otimes id_B}>> \KK \odot B \\
\endCD$$  commutes. Thus, upon completion, it commutes asymptotically.

Property (2) follows because $\{\Phi^{A \otimes B}_t\} = \{\Phi^A_t\}
\otimes \{id_B\}$ as asymptotic morphisms. \qed \enddemo

To finish the proof that $\Phi^A_*$ is an
isomorphism for all $A$, we only need to construct an element $b \in
K_0(C_0(T^*\R^n))$ such that $\Phi^{\C}_*(b) = 1$.

Let $E = \Lambda^*_{\C} \R^n$ denote the complexified exterior algebra of
Euclidean $n$-space. We have that the $C^*$-algebra of endomorphisms of
$E$ is $\End(E) \cong M_{2^n}(\C)$ since $\dim(E) = 2^n$. Let $\epsilon$ be
the grading operator of $E = E^{even} \oplus E^{odd}$ into even and odd forms. For
$v \in \R^n$, let $c(v) : E \to E$ denote the operation defined by
$$c(v)\omega = d_v\omega - \delta_v\omega,\ \omega \in E,$$ where
$d_v = v \wedge$ denotes exterior multiplication by $v$ and $\delta_v = -
d_v^* = v \im $ denotes interior multiplication by $v$. Define an $\R$-linear map 
(also denoted by $c$)
$$c : T^*\R^n \to \End(E)$$ by the formula
$$c(v, \xi) = c(\i \xi) + c(v) = \i (d_\xi + \delta_\xi) + (d_v - \delta_v) ,$$
where we identify $T^*\R^n \cong \R^n \times \R^n$. This ``Clifford
multiplication'' $c$ has the following properties: For all $(v, \xi) \in
T^*\R^n$:
\roster
\item $c(v, \xi)^* = c(v, \xi)$ is self-adjoint.
\item $c(v, \xi)\epsilon = - \epsilon c(v,\xi)$.
\item $c(v,\xi)^2 = \|v\|^2 + \|\xi\|^2$.
\endroster

Our proposed element $b$ will be induced by the following. Define 
the graded $*$-homomorphism 
$$\aligned \Psi : C_0(\R) &\to C_0(T^*\R^n) \otimes \End(E) \\
f &\mapsto f(c), \endaligned$$
where $f(c)(v, \xi) = f(c(v,\xi))$ is defined via the $C^*$-algebra functional calculus.
Since
$$c(v,\xi)^2 = \|v\|^2 + \|\xi\|^2,\quad  (v,\xi) \in T^*\R^n,$$ it follows that
the function
$$(v, \xi) \mapsto f(c(v,\xi))$$
vanishes at infinity on $T^*\R^n$, and so $$\Psi(f) \in C_0(T^*\R^n,
\End(E)) \cong C_0(T^*\R^n) \otimes \End(E)$$ as desired. 

Letting $A$ be trivially graded, $K_0(A)$ is isomorphic to the
group of graded homotopy classes of graded $*$-homomorphisms from $C_0(\R)$ (with the
even/odd grading) to the graded tensor product $A \gtimes M_2(\KK)$ (with
the standard grading) where the addition is given by direct sum (Theorem 4.7 \cite{Trou98}). Thus, 
this graded $*$-homomorphism defines a $K$-theory class
$$b = \lcl \Psi \rcl \in K_0(C_0(T^*\R^n)).$$

\demo{\bf Definition B.3} Let $\{e_j\}_1^n$ be an orthonormal basis for $\R^n$. For
each $t > 0$, define the first order differential operator $B_t : \S(\R^n, E)
\to \S(\R^n, E)$ on the Schwartz space of rapidly decreasing smooth $E$-valued 
functions by 
$$B_t = \sum_{j =1}^n t^{-1}(d_{e_j} + \delta_{e_j})\Dx j + c(v).$$
The operator $B_t$ is formally self-adjoint and we consider it as an
unbounded operator on the Hilbert space  $\H = L^2(\R, E)$ of square-integrable, 
$E$-valued functions. The definition of $B_t$ is independent of the basis $\{e_j\}$.
\enddemo

We collect the facts we will need about the the operator $B_t$ and its
spectral theory in the following theorem \cite{Hig93,Roe88}.

\proclaim{Theorem B.4} $B_t$ is an essentially self-adjoint elliptic operator on
$\H$. Moreover,
\roster
\item $B_t \epsilon = -\epsilon B_t$. Thus, with respect to the
decomposition $E = E^{even} \oplus E^{odd}$, $$B_t = \pmatrix 0 & B_t^- \\
B_t^+ & 0 \endpmatrix.$$
\item The spectrum of $B_t$ is discrete with real eigenvalues
$\lambda_n = \pm\sqrt{2 n t^{-1}}$, where $n \geq 0$.
\item There is an orthonormal basis of $\H$ consisting of eigenfunctions
of $B_t$.
\item $\Ker(B_t)$ is 1-dimensional and spanned by the $0$-form $f(x) = e^{-t\|x\|^2/2}$.
\item The total symbol of $B_t$ is $\sym_{B_t}(v,\xi) = c(v, t^{-1}\xi)$.
\item $B_t$ is $O(n)$-equivariant.
\item $\Index (B_t^+) = +1$.
\endroster
\endproclaim

Using the Spectral Theorem \cite{Gue98}, define a family of $*$-homomorphisms 
$$\{\E^B_t\} : C_0(\R) \to \B(\H) : f \mapsto f(\bar{B}_t).$$ 

\proclaim{Lemma B.5} $\amn {\E^B} {C_0(\R)} {\KK(\H)}$ defines a continuous
family of graded \newline $*$-homomorphisms.
\endproclaim

Extend the index asymptotic morphism $\amn {\Phi} {C_0(T^*\R^n)} {\KK}$
on $\R^n$ to $2^n \times 2^n$ matrices $$\amn {\Phi}
{M_{2^n}(C_0(T^*\R^n))} {M_{2^n}(\KK)}$$ by applying element-wise. Recall
that $\End(E) \cong M_{2^n}(\C)$ under the identification $E =
\Lambda^*_{\C}\R^n \cong \C^{2^n}$.

\proclaim{Proposition B.6} $\Phi^{\C}_*(b) = +1$
\endproclaim

\demo{Proof} Let $\psi \in C_c(\R^n)$. Since $B_t$ is essentially self-adjoint on 
$\H = L^2(\R^{n}, E)$, we have, by an approximation argument similar to the
proof of Theorem  C.10, that
$$\Lim \|\Phi_t(\psi f(c))  - M_\psi f(\bar{B_t})\| = 0$$ for any $f \in
C_0(\R)$ because the total symbol of $B_1$ is $c$. (See also Lemma 9.4 of
\cite{Hig93}.) 

Now, for any $\epsilon > 0$, there is a $\psi \in C_c(\R^n)$ such that $\|(1
- \psi)f(\bar{B_t})\| < \epsilon$. It follows that for all $f \in C_0(\R)$,
$$\Lim \|\Phi_t(f(c)) - f(\bar{B}_t))\| = 0$$ Thus, we have the asymptotic equivalence 
$\Phi_t \circ \Psi \sim_a \E^B_t$. The result now follows from Corollary 4.8 \cite{Trou98} since
$\Index(B_t^+) = +1.$ \qed \enddemo

\proclaim{Theorem B.7 (Bott Periodicity)} $\Phi^A_* : K_0(C_0(T^*\R^n) \otimes A))
\to K_0(A)$ is an isomorphism of abelian groups. \endproclaim

We now turn to the Thom isomorphism.

Let $\pi : E \to X$ be a Hermitian complex vector bundle on the locally
compact topological space $X$. In this section, we will associate to $E
@>{\pi}>> X$ an injective $*$-homomorphism
$$\Psi^E : C_0(\R) \otimes C_0(X) \to C_0(\R) \otimes C_\tau(E),$$
where $C_\tau(E)$ is a $C^*$-algebra Morita equivalent to $C_0(E)$, such
that the induced map on $K$-theory is the {\it topological} Thom
isomorphism $K^0(X) @>{\cong}>> K^0(E)$ \cite{Kar78}.

\demo{\bf Definition B.8} Let $\b E = \pi^*(\Lambda^*E)$ denote the pull-back over $E$
of the exterior algebra bundle $\Lambda^*E$. Let $\epsilon$ denote the
grading of $\b E$ into even and odd forms $$\b E = \b E^{even} \oplus \b
E^{odd} = \pi^*(\Lambda^{even}E) \oplus \pi^*(\Lambda^{odd}E).$$ 
Let $c : E \to \End(\b E)$ denote the canonical section 
$$c(e) = d_e - \delta_e = e \wedge - e \im$$
of the endomorphism bundle $\End(\b E) \to E$, where $d_e$ denotes
exterior multiplication and $\delta_e = -d_e^*$ interior multiplication by
$e$. \enddemo

\proclaim{Lemma B.9} For all $e \in E$, $c(e) : \b E_e \to \b E_e$ satisfies the
following:
\roster
\item $c(e)^* = c(e)$ is self-adjoint.
\item $c(e)^2 = \|e\|^2$.
\item $c(e) \epsilon = - \epsilon c(e)$.
\endroster
\endproclaim

Let $C_\tau(E)$ denote the $C^*$-algebra of bundle endomorphisms
$\alpha : \b E \to \b E$ which vanish at infinity on $E$, under the
pointwise supremum operator norm
$$\|\alpha\| = \sup_{e \in E} \|\alpha(e)\|.$$

\proclaim{Lemma B.10} $C_\tau(E)$ is Morita equivalent to $C_0(E)$. \endproclaim

\proclaim{Corollary B.11} $K_j(C_\tau(E)) \cong K_j(C_0(E)) \cong K^j(E)$ for all $j$.
\endproclaim

\demo{\bf Definition B.12} Define $\Psi^E : C_0(\R) \odot C_0(X) \to C_0(\R)
\otimes C_\tau(E)$ on elementary tensors $f \otimes g$ by the following
$$f \otimes g \mapsto f(\epsilon x +   c)\pi^*(g),$$
where $\pi^*(g) = g \circ \pi$ is the pull-back of $g$ to $E$, and extend
linearly. That is, we have $$\Psi^E  (f \otimes g)(x, e) = f(\epsilon x + 
c(e))g(\pi(e)),$$
for all $x \in \R$ and $e \in E$, where the endomorphism $f(\epsilon x + 
c)$ is defined via the functional calculus. \enddemo

The proof of the following is left to the reader.

\proclaim{Proposition B.13} $\Psi^E$ extends to an injective $*$-homomorphism 
$$\Psi^E : C_0(\R) \otimes C_0(X) \to C_0(\R) \otimes C_\tau(E).$$
\endproclaim

\demo{\bf Example B.14} Suppose $E = X \times \C^n$ is a trivial bundle. Then $$\b
E = X \times \C^n \times \Lambda^*\C^n.$$
Thus, $C_\tau(E) \cong C_0(X) \otimes C_0(\C^n, \End \Lambda^*\C^n)
\cong C_0(X) \otimes M_{2^n}(C_0(\R^{2n}))$ and it follows that
$$\Psi^E  : f \otimes g \mapsto f(\epsilon x +  c(v,\xi)) \otimes g$$
where $(v, \xi) = v + i \xi \in \C^n \cong \R^{2n}$.
\enddemo

The next two lemmas contain the functorial properties of $\Psi^E$ we will
need below.

\proclaim{Lemma B.15} If $F$ is a Hermitian bundle on $X$ isomorphic to $E$, then
the following diagram commutes:
$$\CD
{C_0(\R) \otimes C_0(X)} @>{\Psi^E}>> {C_0(\R) \otimes C_\tau(E)} \\
@VVV                                            @VVV \\
{C_0(\R) \otimes C_0(X)} @>{\Psi^F}>> {C_0(\R) \otimes C_\tau(F)} \\
\endCD$$
\endproclaim

\proclaim{Lemma B.16} If $f : Y \to X$ is a continuous proper map, then the
following diagram commutes:
$$\CD
{C_0(\R) \otimes C_0(X)} @>{\Psi^E}>> {C_0(\R) \otimes C_\tau(E)} \\
@V{f^*}VV                                        @V{f^*}VV \\
{C_0(\R) \otimes C_0(Y)} @>{\Psi^{f^*E}}>> {C_0(\R) \otimes C_\tau(f^*E)}
\\
\endCD$$
\endproclaim

Let $\{B_t\}_{t \in [1,\infty)} : \S(\R^n, \Lambda^*\C^n) \to \S(\R^n,
\Lambda^*\C^n)$ be the family of operators considered in Theorem B.4.
Let $\epsilon$ be the grading operator of $\Lambda^*\C^n$. Using the
Spectral Theorem \cite{Gue98} again, define a family of $*$-homomorphisms
$$\{\Cal A^B_t\} : C_0(\R) \to C_0(\R) \otimes \B(\H) : f \mapsto f(\epsilon x + \bar{B}_t).$$ 
Note that for each $x \in \R$, the operator $\epsilon x + B_t$ is essentially self-adjoint.

\proclaim{Lemma B.17} $\{\Cal A_t^B\}$ is a continuous family of $*$-homomorphisms. \endproclaim

\proclaim{Lemma B.18} Suppose $E = \{pt\}\times \C^n$. Then for all $f \in C_0(\R)$, we
have that $$\Lim \|\Phi_t(f(\epsilon x + c)) - f(\epsilon x + \bar{B}_t)\| =
0.$$
\endproclaim

\demo{Proof} This follows from the discussion in Proposition B.6. \qed \enddemo

\proclaim{Corollary B.19} If $E = \{pt\} \times \C^n$, then the following
diagram asymptotically commutes: 
$$\CD {C_0(\R) \otimes M_{2^n}(C_0(T^*\R^n))} @>{1 \otimes \Phi_t}>>
{C_0(\R) \otimes M_{2^n}(\KK)} \\
@A{\Psi^E}AA  @A{\Cal A^B_t}AA \\
{C_0(\R)} @>{=}>> C_0(\R) 
\endCD$$
\endproclaim

Let $P_t : \H \to \Ker(B_t)$ denote the orthogonal projection onto the
kernel of the operator $B_t$, where $\H = L^2(\R^n, \Lambda^*\C^n)$. By
Theorem 3.7, $\{P_t\}$ is a continuous family of rank one
projections in $\KK(\H)$, the $C^*$-algebra of compact operators on $\H$.

\proclaim{Lemma B.20} The family $\{\Cal A^B_t\}$ is homotopic to the family
$f \mapsto f \otimes P_t$. \endproclaim

\demo{Proof} The homotopy is given by
$$f \mapsto f(\epsilon x + s^{-1}\bar{B}_t),\quad 0 \leq s \leq 1.$$ First,
we need to check continuity in $s$. Considering the factorization
$$\aligned r_\pm(\epsilon x + s^{-1}\bar{B}_t) - &r_\pm(\epsilon x + r^{-
1}\bar{B}_t) = r_\pm(\epsilon x + s^{-1}\bar{B}_t)(r^{-1} - s^{-
1})\bar{B}_t r_\pm(\epsilon x + r^{-1}\bar{B}_t) \\
&= r_\pm(\epsilon x + s^{-1}\bar{B}_t)(1 - rs{-1})(\epsilon x + r^{-
1}\bar{B}_t)r_\pm(\epsilon x + r^{-1}\bar{B}_t) \\
&\quad + r_\pm(\epsilon x + s^{-1}\bar{B}_t)(1 -rs^{-1})(-\epsilon
x)r_\pm(\epsilon x + r^{-1}\bar{B}_t) \endaligned.$$
we see that for each $x \in \R$ $$\|r_\pm(\epsilon x + s^{-1}\bar{B}_t) -
r_\pm(\epsilon x + r^{-1}\bar{B}_t)\| \leq (1-rs^{-1})(1 + |x|).$$
Thus, as $r \to s$, the operator-valued functions $r_\pm(\epsilon x + r^{-
1}\bar{B}_t)$ on the real line converge uniformly to $r_\pm(\epsilon x +
s^{-1}\bar{B}_t)$ on compact subsets. From the inequality
$$\|r_\pm(\epsilon x + s^{-1}\bar{B}_t)\| \leq \frac{1}{\sqrt{x^2 + 1}}$$
these operators are uniformly small on the complement of $[-x ,x]$ for $x
\geq 0$ large. Therefore, they converge uniformly on $\R$.

Since $(\epsilon x + s^{-1}B_t)^2 = x^2 + s^{-2}B_t^2$, it follows that the
eigenvalues of $(\epsilon x + s^{-1}B_t)$ are of the form $$\lambda_n =
\pm\sqrt{x^2 + 2ns^{-2}t^{-1}}$$ by Theorem ~B.4. Thus, as $s \to
0$, the spectrum of $(\epsilon x + s^{-1}B_t)$ corresponding to $n >0$
goes to infinity. This implies that $$f(\epsilon x+ s^{-1}\bar{B}_t) \to f(x)
P_t$$
in norm as $s \to 0$. Note that if $\supp(f) \subset [-a,a]$, then for small
enough $s$, $$f(\epsilon x + s^{-1}\bar{B}_t) = f(x)P_t.$$ The result
follows. \qed \enddemo

Combining the previous two results we obtain the following.

\proclaim{Corollary B.21} The following diagram commutes up to homotopy:
$$\CD {C_0(\R) \otimes M_{2^n}(C_0(T^*\R^n))} @>{1 \otimes \Phi_t}>>
{C_0(\R) \otimes M_{2^n}(\KK)} \\
@A{\Psi^E}AA  @A{1 \otimes P_t}AA \\
{C_0(\R)} @>{=}>> C_0(\R) 
\endCD$$ \endproclaim

Let $E \to X$ be a Hermitian complex vector bundle on $X$. Let $$\Psi^E_* :
K^0(X) \to K^0(E)$$ denote the mapping induced on $K$-theory by the $*$-
homomorphism $$\Psi^E : C_0(\R) \otimes C_0(X) \to C_0(\R) \otimes
C_\tau(E),$$ where we invoke the isomorphisms $K_1(C_0(\R) \otimes
C_0(Y)) \cong K_0(C_0(Y)) \cong K^0(Y)$ for any locally compact space
$Y$.

\proclaim{Theorem B.22 (Thom Isomorphism)} $\Psi^E_* : K^0(X) \to K^0(E)$ is an isomorphism.
\endproclaim

\demo{Proof} Our proof is divided into the following cases.

1.) If $X = \{pt\}$, then $E = \C^n \cong \R^{2n}$, $\b E = \Lambda^*\C^n$
and 
$$\Psi^E : C_0(\R) \to C_0(\R) \otimes M_{2^n}(C_0(\C^n))$$ 
is the map $f \mapsto f(\epsilon x + c)$. Thus, by Corollary B.21, the induced map on $K$-theory is the
Bott Periodicity isomorphism from Theorem B.7.

2.) Suppose $E = X \times \C^n$ is trivial. Example B.14 above then
shows that $\Psi^E = \Psi \otimes id_{C_0(X)}$ where $\Psi$ is the map in
the previous case. The induced map is then seen to be $\Psi^E_* =
\beta_{C_0(X)}$ which is the Bott Periodicity map constructed in Theorem B.1.

3.) If $E$ is trivializable, the result follows from Lemma B.15 and case 2.

4.) Now use the Mayer-Vietoris sequence and the Five Lemma for $X = X_1 \cup X_2$
where $E$ trivializes over $X_1$ and $X_2$.

5.) In general, cover $X$ with open sets $\{X_j\}$ such that over each
$X_j$, $E$ trivializes as $E_j = E|_{X_j} \cong X_j \times \C^n$. An
induction argument using the previous case and the continuity of $K$-theory
$$K^0(X) = \lim \{K^0(X_j) : X_j \subset X \text{ is open}\}$$
finishes the proof. \qed \enddemo

It follows from the proof that $\Psi^E_* : K^0(X) \to K^0(E)$ is, in fact, the
topological Thom Isomorphism \cite{Kar78}.

\centerline{REFERENCES}
\bibliography{MFIT}
\bibliographystyle{plain}
\enddocument
\end